\input ./style/arxiv-general.cfg
\documentclass[aop,MSNbibl,dvips]{arximspdf}
\makeatletter
   \@ifpackageloaded{graphicx}{}{\usepackage{graphicx}}
\makeatother

\doi{10.1214/14-AOP931} 
\volume{43}
\issue{5}
\pubyear{2015}
\firstpage{2205}
\lastpage{2249}

\makeatletter
\newcommand{\rrvert}{\vert}
\newcommand{\llvert}{\vert}
\newproclaim{defn}[teo]{Definition}
\newtheorem{teo}{Theorem}
\newtheorem{prop}[teo]{Proposition}
\newtheorem{lemma}[teo]{Lemma}
\newtheorem{cor}[teo]{Corollary}
\newproclaim{rem}[teo]{Remark}

\newcommand{\Btil}{\widetilde{B}}
\newcommand{\Ytil}{\widetilde{Y}}
\newcommand{\Wtil}{\widetilde{W}}
\renewcommand{\d}{\mathrm{d}}
\newcommand{\dd}{\mathrm{d}}

\newcommand{\cspan}{\overline{\operatorname{span}}}
\newcommand{\Span}{\operatorname{span}}

\def\Xint#1{\mathchoice
{\XXint\displaystyle\textstyle{#1}}%
{\XXint\textstyle\scriptstyle{#1}}%
{\XXint\scriptstyle\scriptscriptstyle{#1}}%
{\XXint\scriptscriptstyle\scriptscriptstyle{#1}}%
\!\int}
\def\XXint#1#2#3{{\setbox0=\hbox{$#1{#2#3}{\int}$ }
\vcenter{\hbox{$#2#3$ }}\kern-.6\wd0}}
\makeatother

\begin{document}
\begin{frontmatter}

\title{Basic properties of critical lognormal multiplicative chaos}
\runtitle{Basic properties of critical lognormal multiplicative chaos\hspace*{12pt}}

\begin{aug}
\author[A]{\fnms{Julien}~\snm{Barral}\ead[label=e1]{barral@math.univ-paris13.fr}},
\author[B]{\fnms{Antti}~\snm{Kupiainen}\corref{}\thanksref{T1,T2}\ead[label=e2]{antti.kupiainen@helsinki.fi}},
\author[B]{\fnms{Miika}~\snm{Nikula}\thanksref{T1}\ead[label=e3]{miika.nikula@helsinki.fi}},
\author[B]{\fnms{Eero}~\snm{Saksman}\thanksref{T1}\ead[label=e4]{eero.saksman@helsinki.fi}}
\and
\author[B]{\fnms{Christian}~\snm{Webb}\thanksref{T1}\ead[label=e5]{christian.webb@helsinki.fi}}
\runauthor{J. Barral et al.} 
\affiliation{Universit\'e Paris 13, University of Helsinki, University
of Helsinki,\\ University of Helsinki and University of Helsinki}
\address[A]{J. Barral\\
LAGA (UMR 7539)\\
D\'epartement de Math\'ematiques\\
Institut Galil\'ee\\
Universit\'e Paris 13\\
99 avenue Jean-Baptiste Cl\'ement\\
93430 Villetaneuse\\
France\\
\printead{e1}}
\address[B]{A. Kupiainen\\
M. Nikula\\
E. Saksman\\
C. Webb\\
Department of Mathematics and Statistics\\
University of Helsinki\\
P.O. Box 68\\
FIN-00014\\
Finland\\
\printead{e2}\\
\phantom{E-mail: }\printead*{e3}\\
\phantom{E-mail: }\printead*{e4}\\
\phantom{E-mail: }\printead*{e5}} 
\end{aug}
\thankstext{T1}{Supported by the Acdemy of Finland.}
\thankstext{T2}{Supported by ERC.}

\received{\smonth{4} \syear{2013}}
\revised{\smonth{3} \syear{2014}}

%
\begin{abstract}
We study one-dimensional exact scaling lognormal multiplicative chaos
measures at criticality. Our main results are the determination of the
exact asymptotics of the right tail of the distribution of the total
mass of the measure, and an almost sure upper bound for the modulus of
continuity of the cumulative distribution function of the measure. We
also find an almost sure lower bound for the increments of the measure
almost everywhere with respect to the measure itself, strong enough to
show that the measure is supported on a set of Hausdorff dimension $0$.
\end{abstract}

%
\begin{keyword}[class=AMS]
\kwd[Primary ]{60G57}
\kwd[; secondary ]{60G18}
\kwd{83C45}
\end{keyword}
\begin{keyword}
\kwd{Multiplicative chaos}
\kwd{critical temperature}
\end{keyword}

\end{frontmatter}

\setcounter{footnote}{2}

\section{Introduction}

Multiplicative chaos is a theory developed by Kahane in the eighties
\cite{K1,K2,K3}. It deals with multiplicative processes generating
martingales, which take values in the cone of nonnegative Radon
measures on $\sigma$-compact metric spaces. This theory is based on the
lognormal multiplicative chaos proposed by Mandelbrot to model
turbulence \cite{M1}, as well as the works previously achieved by
Kahane and Peyri\`ere \cite{kape76} on the simplified model of
multiplicative cascades on trees still proposed by Mandelbrot~\cite
{M2,M3}, namely the so-called Mandelbrot cascades, which assume no
lognormality property. The study of random measures generated by such
multiplicative processes also originates from random covering and
percolation theory questions (see \cite{K4,K2,K3,K5,Fan04,BFan05}).
When statistically self-similar, as it is the case for limits of
Mandelbrot cascades, these measures provide nice illustrations of the
so-called multifractal formalism, as well as models in the study of
intermittent phenomena beyond turbulence, like the distribution of rare
minerals in earth \cite{M89} or stock exchange fluctuations in finance
\cite{M97}. Examples of such measures on $\mathbb{R}^d$ possessing continuous
(rather than only discrete for limits of Mandelbrot cascades) scaling
properties are some of the Gaussian multiplicative chaos built by
Kahane in \cite{K1} or the L\'evy multiplicative chaos built by Fan in
\cite{Fan96}, the compound Poisson cascades built by Barral and
Mandelbrot \cite{BM1} and their generalization to the so-called
infinitely divisible cascades by Bacry and Muzy in \cite{BaMu}.

Kahane's lognormal multiplicative chaos has been recently revisited and
completed in several directions \cite{RoVa1,RoVa2,ARhVa11}. Also, it is
now a central tool in two-dimensional quantum gravity theory since it
provides, through the exponential of the Gaussian free field, the
random measures used to obtain the first rigorous results in direction
to the so-called KPZ formula in works by Duplantier and Sheffield \cite
{PRL,DS08}, as well as Rhodes and Vargas \cite{rhovar08} (see also
Benjamini and Schramm \cite{BS09} for a one-dimensional version in the
framework of Mandelbrot multiplicative cascades on $[0,1]$).
Nondegenerate limits of lognormal multiplicative chaos associated with
the exponential of the Gaussian free field on the circle have also been
used successfully by Astala, Jones, Kupiainen and Saksman in \cite
{AJKS} to build random planar curves by conformal welding. The families
of Gaussian multiplicative chaos considered in these questions are
naturally parameterized  by a continuous parameter $\beta\in[0,\beta_c)$. In the
application to quantum gravity, $\beta$ is in bijection with the
so-called central charge; in random energy models, it corresponds to
the inverse of a temperature; in turbulence, it is a measure of the
intermittence; from a purely geometric viewpoint, it is a decreasing
function of the Hausdorff dimension of the associated measure in the
Euclidean geometry. At the critical temperature, and below it, the
limit $\mu_\beta$ of the martingale $\mu_{\beta,t}$ provided by the
associated multiplicative process vanishes almost surely. For $\beta
>\beta_c$, it is nevertheless possible to give a sense to the
corresponding dual KPZ formula \cite{Duphouches,bajinrhovar12} by
considering measures essentially by subordinating a suitable
nondegenerate Gaussian multiplicative chaos to some stable L\'evy
subordinators; this yields an atomic measure.

At the critical value $\beta_c$, one needs new results in
multiplicative chaos theory. They were recently obtained by Duplantier,
Rhodes, Sheffield and Vargas in \cite{drsv12-1,drsv12-2},
inspired by results recently achieved by A\"id\'ekon and Shi in the
context of the martingales in the branching random walk \cite{aishi11}.
Thus, it is possible to get a nontrivial positive measure at the
critical temperature as the limit of the signed measures $-\frac
{\dd\mu_{\beta,t}}{\dd\beta}|_{\beta=\beta_c}$ as $t\to
\infty$.
Moreover, this measure is continuous. We also mention that like in the
context of martingales in the branching random walks \cite{aishi11},
the critical measure can be obtained as limit in probability of $\mu
_{\beta_c,t}$ properly normalized \cite{drsv12-2}.
During the completion of this paper, corresponding normalization
results \cite{marhva13} were obtained also in the case $\beta> \beta
_c$. These normalization results are analogous to those known in the
branching random walk and random energy models frameworks \cite
{we11,ma11,barhova12}.

This paper is dedicated to the study of some properties of such
critical lognormal multiplicative chaos measure. We concentrate on the
exactly scale-invariant one-dimensional construction. Our main results
are the determination of the asymptotic behavior of the tail of the
distribution of the total mass of the measure, a bound for the modulus
of continuity of the measure for which the previous tail asymptotic
behavior is crucial, and an estimate from below of the measure
increments almost everywhere with respect to the measure, which
completes the estimation provided by the modulus of continuity and goes
beyond the simple fact that the measure has Hausdorff dimension 0; see
Theorems \ref{teotail}, \ref{teomodulus} and \ref
{teohausdorff-gauge} below.

As a motivation to study the exactly scale invariant measure, let us
note that the Gaussian field used to construct the exactly scale
invariant measure in one dimension is simply the Gaussian free field
restricted to a line segment. Thus, the measure can be viewed as a
boundary measure of Liouville quantum gravity (see, e.g., \cite{DS08})
and conjecturally as the boundary measure of random planar maps mapped
to the upper half-plane. Moreover, while the results are mainly stated
for the one-dimensional exactly scale invariant measure, we expect
similar results to hold quite generally for Gaussian multiplicative
chaos measures in any dimension. In Section~\ref{sehigherdimensions},
we finish the paper with a discussion of extensions of our results to a
higher-dimensional setting.


\subsection{Definitions and notation}

In this section, we fix notation and give the precise definitions of
the objects studied in this paper. Formally, the one-dimensional
lognormal multiplicative chaos measures $\mu_\beta$ are random measures
given by
%
%
\begin{equation}
\label{eqformal-exponential} \mu_\beta(\d x) = e^{\beta X(x) - ({\beta^2}/{2}) \mathbb
{E}X(x)^2}\, \d x,
\end{equation}
where $(X(x))_{x \in\mathbb{R}}$ is a logarithmically correlated centered
Gaussian field, that is, a centered Gaussian process with
\[
\mathbb{E}X(x) X(y) \sim\log\frac{1}{|x-y|}\qquad\mbox{as } |x-y| \to0.
\]
However, the logarithmic singularity of the correlation kernel implies
that the realizations of $X$ are not smooth enough to be functions, but
must instead be defined as random distributions. To overcome this major
technical obstacle, in Kahane's theory of multiplicative chaos one
gives a rigorous meaning to the expression (\ref
{eqformal-exponential}) by considering nonsingular approximations $X_t$
to the field $X$, defining the measures $\mu_{\beta,t}$ corresponding
to these regularizations and then taking the weak limit of the measures
$\mu_{\beta,t}$ as the approximation parameter is taken to infinity. In
this way, one completely avoids the problem of defining the exponential
of a distribution.

We mainly concentrate on the \emph{exactly scale invariant}
construction. This scaling property, to be defined below, is central to
the proof of Theorem~\ref{teotail}. The one-dimensional exactly scale
invariant construction is most easily understood through the following
geometric construction, originally due to Bacry and Muzy \cite{BaMu}.

Let $\lambda$ denote the hyperbolic area measure on the upper
half-plane, that is,
\[
\lambda(A) = \int_A \frac{\d x\, \d y}{y^2}\qquad\mbox{for
all } A \subset\mathbb{R}\times\mathbb{R}^+.
\]
For $x \in\mathbb{R}$ and $ t\in\mathbb{R}_+$, let $\mathcal
{C}_t(x)$ denote the set
\[
\mathcal{C}_t(x) = \bigl\{ \bigl(x',y'
\bigr) \vert y' > \max\bigl(2\bigl|x'-x\bigr|, e^{-t}
\bigr), \bigl|x'-x\bigr| < \tfrac{1}{2} \bigr\}
\]
and for a compact interval $I \subset\mathbb{R}$ of length less than
or equal
to 1, denote
\[
\mathcal{C}_t(I) = \bigcap_{x \in I}
\mathcal{C}_t(x).
\]
Note that for $t \geq\log1/|I|$ we have $\mathcal{C}_t(I) = \mathcal
{C}_{\log
1/|I|}(I)$. Next, let $W$ denote the white noise on $\mathbb{R}\times
\mathbb{R}^+$
with control measure $\lambda$. We consider $W$ a random real function
on the Borel sets of $\mathbb{R}\times\mathbb{R}^+$ with finite
$\lambda$-measure
characterized by the following properties: for all disjoint Borel sets
$A, B \subset\mathbb{R}\times\mathbb{R}^+$ such that $\lambda(A),
\lambda(B) < \infty$:
\begin{longlist}[(3)]
\item[(1)] $W(A)$ is a centered Gaussian random variable with variance
$\lambda(A)$,
\item[(2)] the random variables $W(A)$ and $W(B)$ are independent, and
\item[(3)] almost surely we have $W(A \cup B) = W(A) + W(B)$.
\end{longlist}
%
Define
\[
X_t(x) = W\bigl(\mathcal{C}_t(x)\bigr)\qquad\mbox{for
all } x \in\mathbb{R}, t \in[0,\infty).
\]
For a fixed $t > 0$, the covariance structure of the process
$(X_t(x))_{x \in\mathbb{R}}$ can be computed to be
\[
\mathbb{E}X_t(x) X_t(y) = \cases{
\displaystyle t + 1 - e^t |x-y|, &\quad$|x-y| < e^{-t}$,
\vspace*{3pt}\cr
\displaystyle \log\frac{1}{|x-y|}, &\quad$e^{-t} \leq|x-y| \leq1$,
\vspace*{3pt}\cr
0, &\quad$1 <
|x-y|$.}
\]
For any interval $I \subset\mathbb{R}$ of length less than or equal
to $1$ and
$x \in I$, we denote
\[
X_t(I) = W\bigl(\mathcal{C}_t(I)\bigr)\quad\mbox{and}\quad X_t^I(x) =
W\bigl(\mathcal{C}_t(x)
\setminus\mathcal{C}_t(I)\bigr)
\]
to obtain the decomposition $X_t(x) = X_t(I) + X_t^I(x)$, where
$X_t(I)$ is independent of the process $(X_t^I(x))_{x \in I}$. Since
$\mathcal{C}_t(I) = \mathcal{C}_{\log1/|I|}(I)$ for $t \geq\log
1/|I|$, we denote
$X(I):= X_{\log1/|I|}(I)$. Owing to the geometry of the construction,
the field $(X_t(x))$ satisfies the following scale invariance property:
for all intervals $I \subset\mathbb{R}$ and $e^{-t} < |I| < 1$, we have
%
%
\begin{equation}
\label{eqx-exact-scale-invariance} \bigl( X_t(x) \bigr)_{x \in I} =
\bigl(
X_t(I) + X_t^I(x) \bigr)_{x
\in I}
\stackrel{d} {=} \bigl( X_t(I) + X_{t-\log|I|}'
\bigl(x/|I| \bigr) \bigr)_{x \in I},
\end{equation}
where $X'$ is an independent realization of the field $X$. For the
reader's convenience, we give the geometric explanation for this
scaling property in the \hyperref[app]{Appendix}.


For $\beta\in(0,\sqrt{2})$, we construct the measures $\mu_{\beta,t}$
on the unit interval by setting
%
%
\begin{equation}
\label{eqsubcritical-measure} \mu_{\beta,t}(I) = \int_I
e^{\beta X_t(x) - ({\beta^2}/{2})
\mathbb{E}
X_t(x)^2}\, \d x
\end{equation}
for all intervals $I \subset[0,1]$. This construction fits into the
framework of Kahane's theory of multiplicative chaos \cite{K2}, which
implies that almost surely the limit $\mu_\beta= \lim_{t \to\infty}
\mu_{\beta,t}$ exists in the sense of weak convergence of measures and
that the limit measure satisfies $\mu_\beta(I) > 0$ for all intervals
$I \subset[0,1]$. The scaling property (\ref
{eqx-exact-scale-invariance}) implies that the measures $\mu_\beta$ are
exactly scale invariant, especially that
%
%
\begin{equation}
\label{eqsubcritical-exact-scale-invariance}
\qquad \mu_\beta(I) \stackrel{d}
{=}|I| e^{\beta X(I) - ({\beta^2}/{2})
\mathbb{E}X(I)^2} \mu
_\beta'\bigl([0,1]\bigr)\qquad\mbox{for all intervals }
I \subset[0,1],
\end{equation}
where $\mu_\beta'$ is an independent realization of $\mu_\beta$ and
$X(I)$, defined as above, is a centered Gaussian random variable of
variance $\log\frac{1}{|I|}$.

Kahane's work also implies that the corresponding construction for
$\beta\geq\sqrt{2}$ results in degenerate limit measures, that is,
the limit measure will be almost surely null. However, the exact
scaling relation above makes sense for all $\beta> 0$. It has recently
been shown by Duplantier, Rhodes, Sheffield and Vargas \cite{drsv12-1}
that by defining for each interval $I\subset[0,1]$
%
%
\begin{eqnarray}\label{eqcritical-measure}
\mu_{\sqrt{2},t}(I) &=& -\frac{\d}{\d\beta} \bigg|_{\beta=
\sqrt{2}}
\mu_{\beta,t}(I)
\nonumber\\[-8pt]\\[-8pt]
&=& \int_I \bigl(\sqrt{2} (t+1) -X_t(x)\bigr) e^{\sqrt
{2} X_t(x) - \mathbb{E}X_t(x)^2}\, \d x\nonumber
\end{eqnarray}
one\vspace*{2pt} obtains a nondegenerate almost sure weak limit $\mu_{\sqrt{2}} =
\lim_{t \to\infty} \mu_{\sqrt{2},t}$ for which $\mu_{\sqrt
{2}}(I) > 0$
almost surely for all intervals $I \subset[0,1]$. As in the case of
branching random walks (or equivalently, multiplicative cascades), this
derivative turns out to be the correct replacement for the measures
(\ref{eqsubcritical-measure}) in the case $\beta= \sqrt{2}$, at the
very least in the sense that $\mu_{\sqrt{2}}$ is nontrivial and turns
out to satisfy the exact scaling property: as detailed in the \hyperref[app]{Appendix},
we have especially
%
%
\begin{equation}
\label{eqexact-scale-invariance} \qquad\mu_{\sqrt{2}}(I) \stackrel{d} {=}|I|
e^{\sqrt{2} X(I) - \mathbb
{E}X(I)^2}
\mu_{\sqrt
{2}}'\bigl([0,1]\bigr)\qquad\mbox{for all intervals
} I \subset[0,1].
\end{equation}

In defining the lognormal multiplicative chaos measure for the critical
parameter value $\beta= \sqrt{2}$, the peculiar normalizing factor
$(\sqrt{2}(t+1) - X_t(x))$ may also be replaced by a normalization that
is deterministic and also independent of $x$. Inspired by the arguments
of A\"id\'ekon and Shi \cite{aishi11} in the case of branching random
walks, Duplantier, Rhodes, Sheffield and Vargas \cite{drsv12-2}
recently proved that there exists a deterministic constant $c > 0$ such
that for every interval $I \subset[0,1]$ one has
%
%
\begin{equation}
\label{eqrenormalized-critical-measure} \sqrt{t} \int_I e^{\sqrt{2}
X_t(x) - \mathbb{E}X_t(x)^2}\, \d x \to c
\mu_{\sqrt{2}}(I)\qquad\mbox{in probability as } t \to\infty.
\end{equation}

Before moving on to the statements of our results on the fine
properties of~$\mu_{\sqrt{2}}$, we make a final comment on the scale
invariance properties of multiplicative chaos measures. In \cite
{drsv12-1} and \cite{drsv12-2}, the authors deal primarily with a
slightly different construction, the \emph{$\star$-scale invariant}
lognormal multiplicative chaos measures. In terms of the geometric
construction presented here, a $\star$-scale invariant random measure
is obtained by replacing the field $(X_t(x))$ in (\ref
{eqsubcritical-measure}), (\ref{eqcritical-measure}) or (\ref
{eqrenormalized-critical-measure}) by the field $(X_t(x) - X_0(x))$.
Since we will make use of this construction in the proof of Theorem
\ref
{teomodulus}, we have included details on $\star$-scale invariance in
the \hyperref[app]{Appendix}. However, as also noted in the papers themselves, the
proofs of the convergence results in \cite{drsv12-1} and \cite
{drsv12-2} are insensitive to these differences.


\subsection{Main results}

We will make use of the following result of Duplantier, Rhodes,
Sheffield and Vargas \cite{drsv12-2}, which is a corollary of the
deterministic normalization (\ref{eqrenormalized-critical-measure}).

{\renewcommand{\theteo}{\Alph{teo}}
%
%
\begin{teo}\label{teomoments}
For all $h \in(0,1)$, $\mathbb{E}(\mu_{\sqrt{2}}([0,1])^h) < \infty$.
\end{teo}}%

The first of our main theorems is a strengthening of this
result, and analogous to the theorem of Buraczewski \cite{bu09} on the
fixed points of the smoothing transform.

\setcounter{teo}{0}
%
%
\begin{teo}\label{teotail}
The tail probability of $\mu_{\sqrt{2}}$ has the asymptotic behavior
\[
\lim_{\lambda\to\infty}\lambda\mathbb{P} \bigl( \mu_{\sqrt
{2}}
\bigl([0,1]\bigr) > \lambda\bigr) = c_1,
\]
where the constant is given explicitly by
\[
c_1 = \frac{2}{\log2} \mathbb{E}\mu_{\sqrt{2}}\bigl([0,1/2]
\bigr) \log\biggl( 1 + \frac
{\mu_{\sqrt{2}}([1/2,1])}{\mu_{\sqrt{2}}([0,1/2])} \biggr) < \infty.
\]
\end{teo}

This theorem allows one to get detailed information on the geometric
properties of the measure $\mu_{\sqrt{2}}$. The following result is
analogous to our earlier result \cite{bknsw12} on multiplicative cascades.

%
\begin{teo}\label{teomodulus}
For any interval $I\subset[0,1]$ and $\gamma<\frac{1}{2}$,
%
%
\begin{equation}
\mu_{\sqrt{2}}(I)\leq C(\omega) \bigl(\log\bigl(1+|I|^{-1}\bigr)
\bigr)^{-\gamma},
\end{equation}
where $C(\omega) > 0$ is an almost surely finite random constant.
\end{teo}

The proof of this theorem is inspired by the earlier result,
but as the correlations of the field $X$ in the construction of $\mu
_{\sqrt{2}}$ are much more intricate than in the branching random walk
underlying the cascade measures, more involved arguments are needed.

%
\begin{rem}\label{remnoatoms}
We note that this result gives another proof for the result of \cite
{drsv12-1} stating that almost surely, $\mu_{\sqrt{2}}$ has no atoms.
\end{rem}

We also get a bound on the appropriate Hausdorff gauge function to
measure the size of the smallest Borel sets fully supporting $\mu
_{\sqrt
{2}}$. We have the following result.
%

\begin{teo}\label{teohausdorff-gauge}
Denote
$
f_\alpha(n) = \exp( - \sqrt{6 \log2} \sqrt{n (\log n
+\alpha\log\log n )} )
$
for $\alpha> \frac{1}{3}$. Almost surely,
\[
\mu_{\sqrt{2}} \bigl( \bigl\{ x\dvtx  \mu_{\sqrt{2}}\bigl(I_n(x)
\bigr) \geq f_\alpha(n)\mbox{ for all but finitely many }n \bigr\}
\bigr) = \mu_{\sqrt{2}}\bigl([0,1]\bigr),
\]
where $I_n(x) \subset[0,1]$ is the dyadic interval of length $2^{-n}$
containing $x$.
\end{teo}


The proof uses large deviations estimates exploiting both the
exact scaling property of $\mu_{\sqrt{2}}$ and the tail probabilities
given by Theorem~\ref{teotail}. This theorem implies the weaker claim
that almost surely there exists a set of full $\mu_{\sqrt{2}}$-measure
that has Hausdorff dimension $0$, a fact that we state as
Corollary~\ref{cor0-dimension}.

For the log-normal critical Mandelbrot measure $\mu$ on trees, we
establish in \cite{bknsw12} that $
\mu( \{ x\dvtx  \mu(I_n(x)) \leq\psi(n)$ for all but
finitely many $n \} ) = \mu([0,1])
$, for all functions $\psi(n)=n^{-k}$, $k\ge1$. In particular, the
modulus of continuity (shown to be \mbox{optimal}) does not capture the
measure increments behavior $\mu$-almost everywhere---indeed, this is
something one would expect of any multifractal measure. The proof
exploits fine information about the renormalization theory for the low
temperature measures $\mu_{\beta,n}$. Establishing this result in the
present setting remains a challenge, as does proving the optimality of
the bound provided by Theorem~\ref{teomodulus}.


\section{Tail probabilities}

The proof of Theorem \ref{teotail} follows the same idea as the
earlier closely related results of Durrett and Liggett \cite{duli83},
Guivarc'h \cite{gu90}, Liu \cite{li01}, Buraczewski \cite{bu09,bu07}
and Barral and Jin \cite{baji12}: one uses the smoothing transform (or
in the case of multiplicative chaos, a similar distributional equation
with more dependencies) to derive a Poisson equation satisfied by the
quantity one is interested in, and then analyzes the behavior of the
solutions of the Poisson equation at infinity. A key point in the
derivations of the Poisson equations in all these proofs is the use of
an alternate probability measure (the Peyri\`ere probability), the idea
of which goes back to the seminal paper of Kahane and Peyri\`ere \cite
{kape76}. While using quite different kinds of methods, we also point
out the result of Fyodorov and Bouchaud (\cite{fb}), where in the
specific case of the Gaussian free field restricted to the unit circle,
an explicit probability distribution for $\mu_\beta([0,1])$ is obtained (though nonrigorously).

We also note that our form of the tail is related to the freezing
transition scenario: it is believed (see, e.g., \cite{cld}) that a
freezing transition occurs in essentially any logarithmically
correlated random energy model and one universal feature of these
models is that at the critical point, the Laplace transform should be
of the form $1-\mathbb{E}(\exp(-e^{-\beta_c x}\mu_{\beta
_c}([0,1])))\asymp
xe^{-\beta_c x}$ as $x\to\infty$. This is consistent with the tail
$\mathbb{P}(\mu_{\beta_c}([0,1])>x)\asymp x^{-1}$ being universal.

In this section, we denote $\mu:= \mu_{\sqrt{2}}$ and $Y:= \mu
([0,1])$. The variable $Y$ may be written as the fixed point of a
``nonindependent smoothing transform'' as follows:
%
%
\begin{eqnarray}\label{eqsmoothing-transform}
Y &=& \mu\bigl([0,1]\bigr)\nonumber
\\
&=& \mu\bigl([0,1/2]\bigr) + \mu\bigl
([1/2,1]\bigr)
\nonumber\\[-8pt]\\[-8pt]
&=:& \tfrac{1}{2} e^{\sqrt{2} X([0,1/2]) - \mathbb{E}X([0,1/2])^2} Y_0 +
\tfrac{1}{2}
e^{\sqrt{2} X([1/2,1]) - \mathbb{E}X([1/2,1])^2} Y_1
\nonumber
\\
&=:& W_0 Y_0 + W_1 Y_1.\nonumber
\end{eqnarray}
Note that this requires that $\mu(\lbrace\frac{1}{2}\rbrace)=0$ almost
surely. To see this, simply note that by the scaling relation (\ref
{eqexact-scale-invariance}) and Theorem \ref{teomoments} we have, for
any given $h \in(0,1)$, $\mathbb{E}\mu([\frac{1}{2}-\varepsilon,\frac
{1}{2}+\varepsilon])^h \to0$ as $\varepsilon\to0$. By the exact
scale invariance property (\ref{eqexact-scale-invariance}) of $\mu$ we
have $Y_0 \perp W_0$, $Y_1 \perp W_1$ and $Y_0 \stackrel{d}{=}Y
\stackrel{d}{=}Y_1$.
Note, however, that $Y_0$ is not independent of either $Y_1$ or $W_1$.

We then define the version of Peyri\`ere probability that is most
convenient for our needs.
%

\begin{defn}\label{defnpeyriere-probability}
Let $(\Omega,\mathcal{F},\mathbb{P})$ denote the probability space
on which $\mu$
is defined and define a probability space $(\Omega\times\{0,1\},
\mathcal{F}\times\sigma(\{0,1\}), \mathbb{Q})$ by setting
\[
\mathbb{E}_\mathbb{Q}f(\omega,j) = \mathbb{E} \bigl( W_0(
\omega) f(\omega,0) + W_1(\omega) f(\omega,1) \bigr)
\]
for all bounded $\mathcal{F}\times\sigma(\{0,1\})$-measurable
functions $f\dvtx
\Omega\times\{0,1\}$. Define the random variables $\Ytil$, $\Wtil$
and $\Btil$ on this probability space by setting
\[
\Ytil(\omega,j) = Y_j(\omega), \qquad\Wtil(\omega,j) =
W_j(\omega)
\]
and
\[
\Btil(\omega,j) = \cases{
W_1(\omega) Y_1(\omega), &\quad$j=0$,
\vspace*{2pt}\cr
W_0(\omega) Y_0(\omega), &\quad$j=1$.}
\]
\end{defn}

For an intuitive idea of what the measure $\mathbb{Q}$ does, consider the
random probability measure on $[0,1]$ defined by $\frac{\mu(\dd x)}{\mu
([0,1])}$. Then $W_0$ can be seen as the (random) probability that a
point sampled according to this measure is in $[0,\frac{1}{2}]$ and
similarly for $W_1$. So $\mathbb{Q}$ can be seen as a probability distribution
that is obtained by weighting with the information of which half of
$[0,1]$ a point sampled according to $\mu$ is in.

We state the essential properties of the variables defined above as the
following lemma.

%
\begin{lemma}\label{lemmaqb-variables}
The following statements hold:
\begin{longlist}[(3)]
\item[(1)] $\Wtil$ and $\Ytil$ are independent.
\item[(2)] $\Ytil$ (under $\mathbb{Q}$) has the same law as $Y$ (under
$\mathbb{P}$).
\item[(3)] $-\log\Wtil$ is a centered Gaussian with variance $2 \log2$.
\end{longlist}
\end{lemma}

\begin{pf}
Let $f,g\dvtx  \mathbb{R}\mapsto\mathbb{R}$ be bounded and continuous. By direct
computation and the independences $W_0 \perp Y_0$ and $W_1 \perp Y_1$,
\begin{eqnarray*}
\mathbb{E}_\mathbb{Q}f(\Wtil) g(\Ytil) &=&  \mathbb{E} \bigl(
W_0 f(W_0) g(Y_0) + W_1
f(W_1) g(Y_1) \bigr)
\\
&=&  2 \bigl( \mathbb{E}W_0 f(W_0) \bigr) \bigl(
\mathbb{E}g(Y_0) \bigr).
\end{eqnarray*}
By taking $g \equiv1$, we see that $\mathbb{E}_\mathbb{Q}f(\Wtil) =
2 \mathbb{E}W_0
f(W_0)$, and taking $f \equiv1$ yields $\mathbb{E}_\mathbb{Q}g(\Ytil
) = \mathbb{E}g(Y_0) =
\mathbb{E}g(Y)$.
Thus, (2) holds, and moreover,
\[
\mathbb{E}_\mathbb{Q}f(\Wtil) g(\Ytil) = \mathbb{E}_\mathbb
{Q}f(\Wtil) \mathbb{E}_\mathbb{Q}g(\Ytil),
\]
which means that $\Wtil$ and $\Ytil$ are independent as claimed in (1).

The law of $-\log\Wtil$ is easy to identify by computing the moment
generating function. Since $W_0 \stackrel{d}{=}W_1 \stackrel
{d}{=}\frac{1}{2} e^{\sqrt
{2 \log2} N - \log2}$, where $N$ is a standard Gaussian,
\begin{eqnarray*}
\mathbb{E}_\mathbb{Q}e^{t (-\log\Wtil)} &=& \mathbb{E}_\mathbb{Q}
\Wtil^{-t} = \mathbb{E} \bigl( W_0^{1-t} +
W_1^{1-t} \bigr)
\\
&=& 2 \mathbb{E}2^{-(1-t)} e^{(1-t)\sqrt{2 \log2} N - (1-t)\log
2} = 2^{2t
- 1 + (1-t)^2}
\\
&=& e^{t^2 \log2}.
\end{eqnarray*}\upqed
\end{pf}

Define the measure $\nu$ on the positive real axis by setting
%
%
\begin{equation}
\label{eqnu} \int f \,\dd\nu= \mathbb{E}Y f(Y) = \mathbb{E}_\mathbb{Q}
\Ytil f(\Ytil)
\end{equation}
for all continuous functions $f\dvtx  \mathbb{R}^+ \to\mathbb{R}$ with
compact support. The
asymptotics of this measure will be determined through the functions
\[
F_{\alpha,\beta}(x) = \nu\bigl(\bigl(\alpha e^x, \beta e^x\bigr]
\bigr)\qquad\mbox{for } 0 < \alpha< \beta.
\]
In terms of $F_{\alpha,\beta}$, the statement of Theorem \ref{teotail}
is essentially equivalent to the following proposition.

%
\begin{prop}\label{propF-alpha-beta-asymptotics}
Let $F_{\alpha,\beta}$ be defined by $\nu$ as above. Then
\[
\lim_{x\to\infty}F_{\alpha,\beta}(x) = c_1 \log
\frac{\beta
}{\alpha},
\]
where
\[
c_1 = \frac{2}{\log2} \mathbb{E}\mu\bigl([0,1/2]\bigr) \log
\biggl( 1 + \frac
{\mu
([1/2,1])}{\mu([0,1/2])} \biggr) < \infty.
\]
\end{prop}

The first step toward the proof of the proposition above is deriving
the Poisson equation satisfied by $F_{\alpha,\beta}$. Let $\tau$ denote
the law of $-\log\Wtil$. By using the independence of $\Wtil$ and
$\Ytil$, we get
\begin{eqnarray*}
\tau\ast F_{\alpha,\beta}(x) &=&  \int_\mathbb{R}\mathbb
{E}_\mathbb{Q}\Ytil\mathbf{1}_{\{
\Ytil\in(\alpha e^{x+y}, \beta e^{x+y}]\}} \tau(\d y)
\\
&=&  \mathbb{E}_\mathbb{Q}\Ytil\mathbf{1}_{\{\Wtil\Ytil\in
(\alpha e^x, \beta e^x]\}
}
\\
&=&  F_{\alpha,\beta}(x) + \mathbb{E}_\mathbb{Q}\Ytil\mathbf
{1}_{\{\Wtil\Ytil\in
(\alpha e^x, \beta e^x]\}} - \mathbb{E}_\mathbb{Q}\Ytil\mathbf
{1}_{\{\Ytil\in(\alpha
e^x, \beta e^x]\}},
\end{eqnarray*}
where the convolution of the measure $\tau$ with a function $F\dvtx \mathbb
{R}\to\mathbb{R}
$ is defined by
\[
\tau\ast F(x) = \int_\mathbb{R}F(x+y) \tau(\d y) = \int
_\mathbb{R}F(x-y) \tau(\d y).
\]
By using\vspace*{2pt} part (2) of Lemma \ref{lemmaqb-variables}, the distributional
equation (\ref{eqsmoothing-transform}) and the definitions of the
variables $\Wtil$, $\Ytil$ and $\Btil$, the term $\mathbb
{E}_\mathbb{Q}\Ytil\mathbf{1}
_{\{\Ytil\in(\alpha e^x, \beta e^x]\}}$ above may be expressed as
\begin{eqnarray*}
\hspace*{-3pt}&& \mathbb{E}_\mathbb{Q}\Ytil\mathbf{1}_{\{\Ytil\in(\alpha e^x,
\beta e^x]\}}
\\
\hspace*{-3pt}&&\qquad =  \mathbb{E}Y
\mathbf{1}_{\{Y \in(\alpha e^x, \beta e^x]\}}
\\
\hspace*{-3pt}&&\qquad =  \mathbb{E} \bigl( W_0 Y_0 \mathbf{1}_{ \{ \alpha e^x - W_1 Y_1 <
W_0 Y_0 \leq
\beta e^x - W_1 Y_1 \} }
+ W_1 Y_1 \mathbf{1}_{ \{ \alpha e^x - W_0 Y_0 < W_1 Y_1
\leq\beta
e^x - W_0 Y_0 \} } \bigr)
\\
\hspace*{-3pt}&&\qquad = \mathbb{E}_\mathbb{Q}\Ytil\mathbf{1}_{ \{ \alpha e^x -
\Btil< \Wtil\Ytil\leq
\beta e^x - \Btil\} }.
\end{eqnarray*}
The previous computations imply that $F_{\alpha,\beta}$ satisfies the
Poisson equation
%
%
\begin{equation}
\label{eqpoisson} F_{\alpha,\beta}(x) - \tau\ast F_{\alpha,\beta}(x) =
\psi_{\alpha,\beta}(x)
\end{equation}
with the function $\psi_{\alpha,\beta}$ given by
%
%
\begin{equation}
\label{eqpsi-alpha-beta} \psi_{\alpha,\beta}(x) = \mathbb{E}_\mathbb
{Q}\Ytil
\mathbf{1}_{\{ \alpha e^x - \Btil<
\Wtil\Ytil\leq\beta e^x - \Btil\}} - \mathbb{E}_\mathbb{Q}\Ytil
\mathbf{1}_{\{\alpha
e^x < \Wtil\Ytil\leq\beta e^x\}}.
\end{equation}

The desired results on the solutions of this Poisson equation at
infinity could be achieved almost exactly in the same way as in
Buraczewski \cite{bu07}, that is, by building on the general theory
developed by Port and Stone \cite{post69}. The following proposition is
originally due to Buraczewski, but we prefer to give it a
self-contained proof of independent interest that uses only elementary
Fourier analysis.


%
\begin{prop}\label{proppoisson-equation}
Let $\nu$ be a locally finite (nonnegative) Borel measure on
$[0,\infty
)$ that grows at most polynomially in the sense that there exist
$\gamma, C > 0$ such that
\[
\nu\bigl((0,x]\bigr) \leq C(1+x)^\gamma\qquad\mbox{for all } x \geq0.
\]
Define the functions
\[
F_{\alpha,\beta}(x) = \nu\bigl(\bigl(\alpha e^x, \beta e^x\bigr]
\bigr)\qquad\mbox{for all } 0 < \alpha< \beta
\]
and assume that for each $\alpha,\beta$ the function $\psi_{\alpha,\beta
}\dvtx  \mathbb{R}\to\mathbb{R}$ is a bounded and continuous function
indexed by the
parameters $\alpha$ and $\beta$ that satisfies
\[
\int_{-\infty}^\infty\bigl(1+|x|\bigr)^2 \bigl\llvert
\psi_{\alpha,\beta
}(x)\bigr\rrvert\,\d x < \infty
\]
and
\[
\int_{-\infty}^\infty\psi_{\alpha,\beta}(x)\, \d x = 0.
\]
Denote
\[
C_{\alpha,\beta} = \int_{-\infty}^\infty x
\psi_{\alpha,\beta}(x)\, \d x
\]
and assume that the map $(\alpha,\beta) \mapsto C_{\alpha,\beta}$
is continuous.
Finally, let $\tau$ be a Gaussian measure on $\mathbb{R}$, that is,
$\tau$ is
the law of a centered Gaussian random variable with variance $\sigma^2
> 0$.

Then, if $F_{\alpha,\beta}$ satisfies the Poisson equation
\[
F_{\alpha,\beta} - \tau\ast F_{\alpha,\beta} = \psi_{\alpha,\beta},
\]
it has the asymptotics
\[
\lim_{x \to\infty} {F_{\alpha,\beta}(x)} = \frac{2}{\sigma^2}
C_{\alpha,\beta}.
\]
\end{prop}

We split our proof of this proposition into two lemmas and a finalizing
convolution argument.

%
\begin{lemma}\label{lemmapoisson-equation}
Let the function $F\dvtx \mathbb{R}\to\mathbb{R}$ be bounded from below
and satisfy
\[
\lim_{x\to-\infty}F(x)=0.
\]
Assume also that $F$ grows at most exponentially at infinity\footnote
{The assumption of exponential growth is used only to ensure that the
convolution with $\tau$ is well defined.} and solves the Poisson equation
%
%
\begin{equation}
\label{eqlemma-poisson} F - \tau\ast F = \psi,
\end{equation}
where $\tau\sim N(0,\sigma^2)$ is as in Proposition \ref
{proppoisson-equation}
and $\psi\dvtx \mathbb{R}\to\mathbb{R}$ satisfies
\[
\int_{-\infty}^\infty\bigl(1+|x|\bigr)^2\bigl|\psi(x)\bigr|\, \d
x <\infty, \qquad\int_{-\infty}^\infty\psi(x) \,\d x = 0
\quad\mbox{and}\quad\lim_{x\to\pm\infty}\psi(x)=0.
\]
Then $F$ has the asymptotics
%
%
\begin{equation}
\label{eqpoisson-lemma-asymptotics} \lim_{x \to\infty} F(x) = \frac
{2}{\sigma^2} \int
_{-\infty
}^\infty x \psi(x) \,\d x.
\end{equation}
%
\end{lemma}

\begin{pf}
To shorten the notation, we denote $\int_{-\infty}^\infty x \psi(x)
\,\d
x = A$.

We start by proving that equation (\ref{eqlemma-poisson}) has some
bounded solution $F_1$ that has the desired asymptotics
%
%
\begin{equation}
\label{eqfirst} \lim_{x\to-\infty} F_1(x) = 0 \quad\mbox{and}\quad\lim
_{x \to
\infty} F_1(x) =
\frac{2}{\sigma^2} A.
\end{equation}
We first consider the case $A=0$. Then our assumptions imply that the
Fourier transform of $\psi$ satisfies [our convention for the Fourier
transform of $\psi\in L^1(\mathbb{R})$ is $\widehat{\psi}(\xi
)=\int_\mathbb{R}
e^{-ix\xi}\psi(x) \,\d x$]
\[
\widehat\psi\in C^2(\mathbb{R}) \cap L^\infty(\mathbb{R})
\quad\mbox{and} \quad\widehat\psi(0) = \widehat\psi'(0) = 0.
\]
As $1-\widehat\tau(\xi)=1-\exp(-\sigma^2 \xi^2/2)$ is smooth with zero
of order 2 at the origin we may directly define $F_1$ in the
distribution sense through
\[
\widehat F_1(\xi):= \frac{\widehat\psi(\xi)}{1-e^{-\sigma^2 \xi
^2/2}} = \widehat\psi(\xi) +
\biggl(\frac{\widehat\psi(\xi)}{\xi^2} \biggr) \biggl(\frac{\xi^2
e^{-\sigma^2 \xi^2/2}}{1-e^{-\sigma^2 \xi^2/2}} \biggr) =: \widehat
\psi(\xi) + \widehat F_2(\xi).
\]
Since obviously $\widehat F_2\in L^1$, we have $\lim_{x \to\pm\infty}
F_2(x) = 0$ by the Riemann--Lebesgue lemma, and the same follows for
$F_1$ by the assumption on $\psi$.

In order to construct a solution $F_1$ in the general case $A\neq0$,
first define $\psi_0 = F_0-\tau*F_0$ with $F_0 = \chi_{(0,\infty)}$,
with $\chi$ referring to the indicator function of a set. Directly from
the definition we see that $\psi_0 \in L^\infty(\mathbb{R})$ and
that $\psi_0$
decays exponentially as $x\to\pm\infty$, so that it satisfies the
moment conditions of the lemma, and moreover $\widehat\psi_0 \in
C^\infty$. Also, $\widehat F_0 (\xi) = \pi\delta_0-i\xi^{-1}$
(here $\xi
^{-1}$ is understood as a principal value distribution). Since
$1-\widehat\tau(\xi)=\sigma^2\xi^2/2+O (\xi^3 )$, we
see that
\[
\widehat\psi_0(\xi) = - \frac{i \sigma^2}{2} \xi+ O\bigl(
\xi^2\bigr)
\]
at the origin. Hence,
\[
\int_{-\infty}^\infty\psi_0(x) \,\d x =
\widehat\psi_0(0) = 0 \quad\mbox{and} \quad\int
_{-\infty}^\infty x \psi_0(x) \,\d x = i
\widehat\psi_0'(0) = \frac{\sigma^2}{2}.
\]

Thus, in the case $A \neq0$ we define $F_1$ by finding the solution
for the Poisson equation (\ref{eqlemma-poisson}) with the right-hand
side $\widetilde\psi= \psi- \frac{2}{\sigma^2} A \psi_0$ and then\vspace*{1pt}
adding $\frac{2}{\sigma^2} A F_0$. At this point, it is clear that the
solution obtained this way is bounded and has the desired behavior at
$\pm\infty$.

Let us finally assume that $F$ and $\psi$ are as in the theorem. Let
$F_1$ be the bounded solution of (\ref{eqlemma-poisson}) constructed
above, so that $F_1$ satisfies the conclusion of the theorem. It is
enough to verify that $H:=F-F_1$ is constant since then $H\equiv0$ by
considering the limit at $-\infty$. Now $H$ is bounded from below and
satisfies the homogeneous Poisson equation
%
%
\begin{equation}
\label{eqhomogeneous-poisson} H = \tau\ast H.
\end{equation}
The claim follows from Lemma \ref{lemmahomogeneous-poisson} below.
\end{pf}

%
\begin{lemma}\label{lemmahomogeneous-poisson}
Let $H$ solve the homogeneous Poisson equation (\ref
{eqhomogeneous-poisson}) and assume that it is bounded from below and
has at most exponential growth at $\pm\infty$. Then $H$ is constant.
\end{lemma}

\begin{pf}
By adding a constant, we may without loss of generality assume that $H
\geq0$. Let $u(x,t)$ ($x\in\mathbb{R}$, $t\geq0$) denote the heat
extension of
$H$ to the upper half-plane, explicitly given by
\[
u(x,t) = \frac{1}{\sqrt{2 \pi t}} \int_{-\infty}^\infty
e^{{-(y-x)^2}/{(2t)}} H(y) \,\d y.
\]
By assumption, $u$ is periodic in $t$:  denoting $t_0 = \sigma^2$,
$u(x,t+t_0)=u(x,t)$ for all $t \geq0$. Define the function $v$ in the
upper-half plane by setting
\[
v(x,t):=\int_0^{t_0}u(x,t+s) \,\d s.
\]
Then $v$ solves the heat equation and, by the periodicity of $u$, it is
constant in $t$. Thus, it is harmonic in $x$, that is, a linear
function $v(x,t)=ax+b$. Here, $a=0$ by the nonnegativity of $u$, whence
$v$ is constant. This shows that there is a constant $C$ independent of
$x$ so that $\int_{t_0/2}^{t_0}u(x,s)\,\dd s <C$. Especially,\vspace*{-2pt} for each $x$
there is $t_1=t_1(x)\in(t_0/2,t_0)$ so that $u(x,t_1)\leq2C/t_0$.
The heat kernel $(2 \pi t)^{-1/2} e^{-x^2/2t}$ can be bounded from
below on $x \in[-1,1]$ uniformly in $t \in(t_0/2,t_0)$,
whence again using the nonnegativity of $H$ we deduce that $\int
_{x-1}^{x+1} H(y) \,\d y \leq C'$ for all $x\in\mathbb{R}$. As we
combine this
information with the fact that $H = \tau\ast H$ it follows that $H$ is
bounded. Then the equation
\[
\bigl(1-e^{-\xi^2/2}\bigr)\widehat H(\xi) = 0,
\]
interpreted in the sense of distributions, shows that $\widehat
H=c_1\delta_0+c_2\delta'_0$,
that is, $H$ is linear. By nonnegativity, we finally deduce that $H$ is
constant.
\end{pf}

%
\begin{rem}\label{rempoissonalternative}
As pointed out by one of the referees, these types of results often
have more probabilistic proofs as well. For example, let us sketch one
for the previous lemma. Consider again $H\geq0$ and note that the
condition $H=\tau*H$ means that $(H(S_n))_{n \geq0}$ is a martingale
for the Gaussian random walk $(S_n)$ with increments distributed
according to $\tau$. Since $H$ is nonnegative, the martingale converges
to some nonnegative random variable, say $\mathcal{H}$. On the other
hand, $(S_n)$ is neighborhood recurrent, so for any $\epsilon>0$ and
$x\in\mathbb{R}$ we can find a subsequence $n_k$ such that $n_k\to
\infty$ as
$k\to\infty$ and $|S_{n_k}-x|<\epsilon$ for all $k$. Now, the fact that
$H = \tau\ast H$ together with the growth condition assumed of $H$
implies that $H$ is a smooth function. Thus, for any given $x$ we have
$\mathcal{H} = \lim_{k \to\infty} H(S_{n_k}) = H(x)$ and, therefore,
$H$ is constant.
\end{rem}

We finish the proof of Proposition \ref{proppoisson-equation} by
deducing it from Lemma \ref{lemmapoisson-equation} by a convolution
argument analogous to the one of Buraczewski \cite{bu07}.

\begin{pf*}{Proof of Proposition \ref{proppoisson-equation}}
Let $\phi\geq0$ be an arbitrary symmetric smooth test function with
$\operatorname{supp} \phi\subset[-1,1]$ and $\int_R \phi= 1$. Given any
locally integrable \mbox{$g\dvtx \mathbb{R}\to\mathbb{R}$}, let $g_\varepsilon$
denote the
convolution $g_\varepsilon= g \ast\varepsilon^{-1}\phi(\varepsilon
^{-1}\cdot)$. By convolving the Poisson equation, we obtain [writing,
e.g., $(F_{\alpha,\beta})_\varepsilon= F_{\alpha,\beta,\varepsilon}$]
for any $0 \leq\alpha< \beta$ and $\varepsilon> 0$
\[
F_{\alpha,\beta,\varepsilon} = \tau\ast F_{\alpha,\beta,\varepsilon} +
\psi_{\alpha,\beta,\varepsilon}.
\]
By the continuity of $\psi_{\alpha,\beta,\varepsilon}$ and
integrability of $\psi_{\alpha,\beta}$, we have $\lim_{x\to\pm
\infty}
\psi_{\alpha,\beta,\varepsilon} = 0$. From Lemma \ref
{lemmapoisson-equation}, we thus obtain, for each $\varepsilon>0$,
the asymptotics
%
%
\begin{equation}
\label{eqpoisson-epsilon-asymptotics} \lim_{x\to\pm\infty} F_{\alpha,\beta,\varepsilon}(x) =
\frac
{2}{\sigma
^2} C_{\alpha,\beta},
\end{equation}
since
\[
\int_\mathbb{R}x \psi_{\alpha,\beta,\varepsilon}(x)\,\d x = i \widehat\psi
_{\alpha,\beta,\varepsilon}'(0) = i \widehat\psi_{\alpha,\beta
}'(0)
= \int_\mathbb{R}x \psi_{\alpha,\beta}(x) = C_{\alpha,\beta}.
\]

In order to remove the $\varepsilon$ from (\ref
{eqpoisson-epsilon-asymptotics}), let $k \in(1, (\beta/\alpha)^{1/2})$
be given and observe that by the definition of $F_{\alpha,\beta}$ as a
measure\vspace*{1pt} of an interval we have $F_{k\alpha,k^{-1}\beta}(x)\leq
F_{\alpha,\beta,\varepsilon}(x)$ for all $x$ as soon as $\varepsilon$
is small
enough. Hence, we obtain from (\ref{eqpoisson-epsilon-asymptotics}) that
\[
\limsup_{x\to\infty} F_{\alpha,\beta}(x) \leq\frac{2}{\sigma^2}
C_{k^{-1} \alpha,k \beta}.
\]
By letting $k\to1^+$ and recalling the assumption of the continuity of
$(\alpha,\beta) \mapsto C_{\alpha,\beta}$ it follows that $\limsup_{x\to
\infty} F_{\alpha,\beta}(x) \leq2C_{\alpha,\beta}/\sigma^2$. The
converse direction $\liminf_{x\to\infty} F_{\alpha,\beta}(x) \geq
2C_{\alpha,\beta}/\sigma$ is obtained analogously by starting from the
inequality $F_{k^{-1}\alpha,k\beta}(x) \leq F_{\alpha,\beta,\varepsilon}(x)$.
\end{pf*}

The proof of Proposition \ref{propF-alpha-beta-asymptotics} has now
essentially been reduced to checking that the Poisson equation (\ref
{eqpoisson}) with $F_{\alpha,\beta}$ determined by $\nu$ as in
(\ref{eqnu}) and $\psi_{\alpha,\beta}$ given by~(\ref{eqpsi-alpha-beta})
satisfies the assumptions of Proposition \ref{proppoisson-equation}.

\begin{pf*}{Proof of Proposition \ref{propF-alpha-beta-asymptotics}}
We first check that the measure $\nu$ satisfies $\nu((0,x]) \leq
C(1+x)^\gamma$ for some $C, \gamma> 0$. This is clear from the
definition and Theorem \ref{teomoments}: for any $\gamma\in(0,1)$,
\[
\nu\bigl((0,x]\bigr) = \mathbb{E}Y \mathbf{1}_{\{Y \in(0,x]\}} \leq
x^\gamma
\mathbb{E}Y^{1-\gamma} \mathbf{1}_{\{Y \in(0,x]\}} \leq\mathbb
{E}Y^{1-\gamma}
x^\gamma.
\]

To check the integrability conditions on $\psi_{\alpha,\beta}$, we
define the functions
\[
\psi_\alpha(x) = \mathbb{E}_\mathbb{Q}\Ytil
\mathbf{1}_{
\{ \alpha e^x - \Btil<
\Wtil\Ytil\leq\alpha e^x \} } \quad\mbox{and}\quad\psi_\beta(x) =
\mathbb{E}_\mathbb{Q}\Ytil\mathbf{1}_{ \{
\beta e^x - \Btil< \Wtil
\Ytil\leq\beta e^x \} }.
\]
By this definition,
\[
\psi_{\alpha,\beta}(x) = \mathbb{E}_\mathbb{Q}\Ytil\mathbf
{1}_{\{ \alpha e^x - \Btil<
\Wtil\Ytil\leq\beta e^x - \Btil\}} - \mathbb{E}_\mathbb{Q}\Ytil\mathbf
{1}_{\Wtil
\Ytil\in(\alpha e^x, \beta e^x]}
= \psi_{\alpha}(x) - \psi_{\beta}(x).
\]
Since the functions $\psi_{\alpha}$ and $\psi_{\beta}$ are
positive, to
check the integrability conditions of Theorem \ref
{proppoisson-equation} on the functions $\psi_{\alpha,\beta}$ it is
sufficient to show that
\[
\int_{-1}^1 \psi_{\alpha}(x) \,\d x <
\infty\quad\mbox{and} \quad\int_{-\infty}^\infty
x^2 \psi_\alpha(x) \,\d x < \infty\qquad\mbox{for all }
\alpha> 0.
\]
In our situation, the first condition is clear, since $\mathbb
{E}_\mathbb{Q}\Wtil
^{-1} < \infty$. For the second condition, some computation and a
separate lemma is needed. We write
\[
\mathbf{1}_{ \{ \alpha e^x - \Btil< \Wtil\Ytil\leq\alpha e^x
\} } = \mathbf{1}_{ \{ ({\Wtil\Ytil})/{\alpha} \leq t <
({\Wtil
\Ytil+ \Btil})/{\alpha} \}} \qquad\mbox{for } t =
e^x
\]
and use the integral
\[
\int_a^b \frac{\log t}{t} \,\d t =
\frac{1}{2} \biggl(\log\frac
{b}{a} \biggr) (\log ab )\qquad\mbox{for } 0 < a < b
\]
to compute
\begin{eqnarray*}
&& \int_0^\infty x^2
\psi_\alpha(x) \,\d x
\\
&&\qquad = \mathbb{E}_\mathbb{Q}\Ytil\int_0^\infty
x^2 \mathbf{1}_{
\{ \alpha e^x - \Btil<
\Wtil\Ytil\leq\alpha e^x \} } \,\d x
\\
&&\qquad = \mathbb{E}_\mathbb{Q}\Ytil\int_1^\infty
\mathbf{1}_{ \{
({\Wtil\Ytil})/{\alpha}
\leq t < ({\Wtil\Ytil+ \Btil})/{\alpha} \}} \frac{(\log
t)^2}{t} \,\d t
\\
&&\qquad = \mathbb{E}_\mathbb{Q}\Ytil\mathbf{1}_{ \{ ({\Wtil
\Ytil})/{\alpha} > 1
\} } \int
_{({\Wtil\Ytil})/{\alpha}}^{({\Wtil\Ytil+ \Btil
})/{\alpha}} \frac{(\log t)^2}{t} \,\d t
\\
&&\quad\qquad{} + \mathbb{E}_\mathbb{Q}\Ytil\mathbf{1}_{
\{ ({\Wtil\Ytil
})/{\alpha} < 1 < ({\Wtil\Ytil+ \Btil})/{\alpha} \} }
\int_1^{({\Wtil\Ytil+ \Btil})/{\alpha}} \frac{(\log t)^2}{t} \,\d t
\\
&&\qquad \leq\frac{1}{2} \mathbb{E}_\mathbb{Q}\Ytil\mathbf{1}_{
\{({\Wtil\Ytil
})/{\alpha} > 1 \}}
\log\biggl( 1 + \frac{\Btil}{\Wtil\Ytil} \biggr) \log\biggl( \frac
{\Wtil\Ytil}{\alpha} \cdot
\frac{\Wtil\Ytil+
\Btil
}{\alpha} \biggr) \log\biggl( \frac{\Wtil\Ytil+ \Btil}{\alpha} \biggr)
\\
&&\quad\qquad{} + \frac{1}{2} \mathbb{E}_\mathbb{Q}\Ytil
\mathbf{1}_{ \{({\Wtil
\Ytil})/{\alpha} < 1 < ({\Wtil\Ytil+ \Btil})/{\alpha} \}} \biggl( \log
\biggl( \frac{\Wtil\Ytil+ \Btil}{\alpha} \biggr)
\biggr)^2 \log\biggl( \frac{\Wtil\Ytil+ \Btil}{\alpha} \biggr)
\\
&&\qquad =: I_1 + I_2,
\end{eqnarray*}
and similarly by the change of variables $s = e^{-x}$ we get
\begin{eqnarray*}
&& \int_{-\infty}^0  x^2
\psi_\alpha(x) \,\d x
\\
&&\qquad = \mathbb{E}_\mathbb{Q}\Ytil\int_{-\infty}^0
x^2 \mathbf{1}_{
\{ \alpha e^x - \Btil
< \Wtil\Ytil\leq\alpha e^x \} } \,\d x
\\
&&\qquad = \mathbb{E}_\mathbb{Q}\Ytil\int_1^\infty
\mathbf{1}_{ \{
\alpha/s - \Btil< \Wtil
\Ytil\leq\alpha/s \} } \frac{(\log s)^2}{s} \,\d s
\\
&&\qquad \leq\frac{1}{2} \mathbb{E}_\mathbb{Q}\Ytil\mathbf{1}_{
\{ ({\Wtil\Ytil+
\Btil})/{\alpha} < 1 \}}
\log\biggl( 1 + \frac{\Btil}{\Wtil
\Ytil} \biggr) \log\biggl( \frac{\alpha}{\Wtil\Ytil} \cdot
\frac{\alpha
}{\Wtil\Ytil+ \Btil} \biggr) \log\biggl( \frac{\alpha}{\Wtil
\Ytil} \biggr)
\\
&&\quad\qquad{} + \frac{1}{2} \mathbb{E}_\mathbb{Q}\Ytil
\mathbf{1}_{ \{({\Wtil
\Ytil})/{\alpha} < 1 < ({\Wtil\Ytil+ \Btil})/{\alpha} \}} \biggl( \log
\biggl( \frac{\Wtil\Ytil}{\alpha} \biggr)
\biggr)^2 \log\biggl( \frac{\alpha}{\Wtil\Ytil} \biggr)
\\
&&\qquad =: I_3 + I_4.
\end{eqnarray*}
To show that $I_1 < \infty$, we use the crude estimate $\log(1+x)
\leq
C_p x^p$, valid for all $p > 0$ for sufficiently large constant $C_p >
0$ depending only on $p$, to get
\begin{eqnarray*}
I_1 &\leq& \frac{C_{p_1} C_{p_2} C_{p_3}}{\alpha^{p_1} \alpha^{p_2}}
\mathbb{E} \biggl( \mu\bigl([0,1/2]\bigr) \biggl( \frac{\mu
([1/2,1])}{\mu
([0,1/2])}
\biggr)^{p_1}
\\
&&\hspace*{65pt}{}\times  \bigl( \mu\bigl([0,1/2]\bigr)^{p_2} + \mu
\bigl([0,1]\bigr)^{p_2} \bigr) \mu\bigl([0,1]\bigr)^{p_3}
\biggr).
\end{eqnarray*}
In Lemma \ref{lemmamixedmoments} below, we show that for any $0 < h <
1$ we have
\[
\mathbb{E}\mu\bigl([0,1/2]\bigr)^h \mu\bigl([1/2,1]
\bigr)^h < \infty.
\]
By choosing $p_1, p_2, p_3 > 0$ such that $0 < 1-p_1+p_2+p_3 < 1$ and
$p_1+p_2+p_3 < 1$, this implies the finiteness of $I_1$. For $I_2$ one
may estimate $\frac{\Wtil\Ytil+ \Btil}{\alpha} \leq1 + \frac
{\Btil
}{\Wtil\Ytil}$ and proceed as in the case of $I_1$. In estimating
$I_3$, one may write $\frac{\alpha}{\Wtil\Ytil+ \Btil} < \frac
{\alpha
}{\Wtil\Ytil}$ and proceed as before, and the finiteness of $I_4$
follows the same route.


In order to apply Proposition \ref{proppoisson-equation}, we still
need to show that
\[
\int_{-\infty}^\infty\psi_{\alpha,\beta}(x) \,\d x = 0
\]
and compute the value of the integral
\[
C_{\alpha,\beta} = \int_{-\infty}^\infty x
\psi_{\alpha,\beta
}(x) \,\d x.
\]
The first integral follows immediately from the integrability of $\psi
_\alpha$ and the fact that
\[
\psi_{\alpha,\beta}(x) = \psi_\alpha(x) - \psi_\beta(x) =
\psi_\alpha(x) - \psi_{\alpha}\biggl(x + \log\frac{\alpha}{\beta}
\biggr).
\]
%
The value of $C_{\alpha,\beta}$ can be calculated by using the change
of variables $x = e^t$ as above to obtain
\[
\int_{-\infty}^\infty x \psi_\alpha(x) \,\d x =
\frac{1}{2} \mathbb{E}_\mathbb{Q}\Ytil\log\biggl( 1 +
\frac{\Btil}{\Wtil\Ytil} \biggr) \log\frac
{\Wtil
\Ytil( \Wtil\Ytil+ \Btil)}{\alpha^2},
\]
which implies
\[
C_{\alpha,\beta} = \int_{-\infty}^\infty x \bigl(
\psi_\alpha(x) - \psi_\beta(x)\bigr) \,\d x =
\mathbb{E}_\mathbb{Q}\Ytil\log\biggl( 1 + \frac{\Btil
}{\Wtil\Ytil} \biggr) \log
\frac{\beta}{\alpha}.
\]
Proposition \ref{proppoisson-equation} now gives the desired asymptotics
\begin{eqnarray*}
F_{\alpha,\beta}(x) & \stackrel{x \to\infty} {\longrightarrow} &
\frac
{2 \mathbb{E}_\mathbb{Q}\Ytil\log( 1 + {\Btil
}/({\Wtil\Ytil}) )}{2
\log2} \log\frac{\beta}{\alpha}
\\
&=&  \frac{2}{\log2} \mathbb{E}\mu\bigl([0,1/2]\bigr) \log\biggl( 1 +
\frac
{\mu
([1/2,1])}{\mu([0,1/2])} \biggr) \log\frac{\beta}{\alpha}
\end{eqnarray*}
for all $0 < \alpha< \beta$.
\end{pf*}

Before moving on to the final step of the proof of Theorem \ref
{teotail}, we complete the proof of Proposition \ref
{propF-alpha-beta-asymptotics} by proving Lemma \ref{lemmamixedmoments}.

%
\begin{lemma}\label{lemmamixedmoments-disjoint}
For any $h \in(0,1)$ and any pair of intervals $I_1, I_2 \subset
[0,1]$ such that $d(I_1,I_2) > 0$,
\[
\mathbb{E} \bigl( \mu(I_1)^h \mu(I_2)^h
\bigr) < \infty.
\]
\end{lemma}

\begin{pf}
We use Kahane's convexity inequality, to be given as Proposition \ref
{propkahane-inequality}, and the definition (\ref
{eqrenormalized-critical-measure}) of the critical measure as the
limit of
\[
\mu_t(\d x) = \sqrt{t} e^{\sqrt{2} X_t(x) - \mathbb{E}X_t(x)^2} \,\d x
\]
as $t \to\infty$. Write the product $\mu_t(I_1) \mu_t(I_2)$ as
\[
\mu_t(I_1) \mu_t(I_2) = t \int
_{I_1} \d x \int_{I_2} \d y\,
e^{\sqrt{2}
(X_t(x) + X_t(y)) - \mathbb{E}X_t(x)^2 - \mathbb{E}X_t(y)^2}
\]
and consider the Gaussian fields $Z_t(x,y) = X_t(x) + X_t(y)$ and
$\widetilde{Z}_t(x,y) = X_t(x) + \widetilde{X}_t(y)$ indexed by $I_1
\times I_2$, where $\widetilde{X}_t$ is an independent realization of
the field $X_t$. The covariance kernel of $\widetilde{Z}_t$ is clearly
dominated by the covariance kernel of $Z_t$, so Proposition \ref
{propkahane-inequality} gives the inequality
\begin{eqnarray*}
&& \mathbb{E} \biggl( t \int_{I_1} \d x \int
_{I_2} \d y \, e^{\sqrt
{2} Z_t(x,y) -
\mathbb{E}Z_t(x,y)^2} \biggr)^h
\\
&&\qquad \leq
\mathbb{E} \biggl( t \int_{I_1} \d x \int_{I_2}
\d y\, e^{\sqrt{2} \widetilde{Z}_t(x,y) - \mathbb{E}\widetilde
{Z}_t(x,y)^2} \biggr)^h
\\
&&\qquad = \mathbb{E} \biggl( \sqrt{t} \int_{I_1} e^{\sqrt{2} X_t(x) -
\mathbb{E}X_t(x)^2}
\,\d x \biggr)^h \mathbb{E} \biggl( \sqrt{t} \int_{I_2}
e^{\sqrt{2}
\widetilde
{X}_t(x) - \mathbb{E}\widetilde{X}_t(x)^2} \,\d x \biggr)^h
\\
&&\qquad = \mathbb{E}\bigl(\mu_t(I_1)^h\bigr)
\mathbb{E}\bigl(\mu_t(I_2)^h\bigr)
\\
&&\qquad  < \infty.
\end{eqnarray*}
By expanding the variance $\mathbb{E}Z_t(x,y)^2 = \mathbb{E}X_t(x)^2
+ \mathbb{E}X_t(y)^2 +
2 \mathbb{E}X_t(x) X_t(y)$ we note that the first expression may be estimated
from below by
\[
e^{-2 \sup_{x \in I_1, y \in I_2} \mathbb{E}X_t(x) X_t(y)} \mathbb
{E}\bigl(t \mu_t(I_1) \mu
_t(I_2)\bigr)^h.
\]
Since the intervals $I_1$ and $I_2$ are separated by a positive
distance, the supremum in the exponent stays bounded as $t \to\infty$,
which proves the claim.
\end{pf}


%
\begin{lemma}\label{lemmamixedmoments}
For any $h \in(0,1)$,
\[
\mathbb{E} \bigl( \mu\bigl([0,1/2]\bigr) \mu\bigl([1/2,1]\bigr)
\bigr)^h < \infty.
\]
\end{lemma}

\begin{pf}
Fix $h \in(0,1)$. For every $k \in\mathbb{N}$, let $J_k = [
1/2 -
2^{-k}, 1/2 + 2^{-k} ]$. Denote the left and right half of $J_k$
by $J_k^0$ and $J_k^1$, and the right and left halves of $J_k^0$ (and
$J_k^1$) by $J_k^{00}$ and $J_k^{01}$ ($J_k^{10}$ and $J_k^{11}$).
Define the sets $A_k$ by
\[
A_k = \bigl(J_k^{00} \times
J_k^{11}\bigr) \cup\bigl(J_k^{00}
\times J_k^{10}\bigr) \cup\bigl(J_k^{01}
\times J_k^{11}\bigr).
\]
Write
\[
Z = \mu\bigl([0,1/2]\bigr) \mu\bigl([1/2,1]\bigr) = \int_{[0,1/2]}
\mu(\d x) \int_{[1/2,1]} \mu(\d y)
\]
and define the random variables
\begin{eqnarray*}
Z_k &=&  \int_{[0,1/2]} \mu(\d x) \int
_{[1/2,1]} \mu(\d y) \chi_{A_k}(x,y)
\\
&=&  \mu\bigl(J_k^{00}\bigr) \mu\bigl(J_k^{11}
\bigr) + \mu\bigl(J_k^{00}\bigr) \mu\bigl(J_k^{10}
\bigr) + \mu\bigl(J_k^{01}\bigr) \mu\bigl(J_k^{11 }
\bigr)
\end{eqnarray*}
for $k \in\mathbb{N}$. It is clear that
\[
Z = \sum_{k=1}^\infty Z_k
\]
and thus by the subadditivity of $x \mapsto x^h$
\[
\mathbb{E}Z^h \leq\sum_{k=1}^\infty
\mathbb{E}Z_k^h.
\]
By the exact scaling property of the construction, the measure $\mu$ satisfies
%
%
\begin{equation}
\label{eqpiecewise-exact-scale-invariance} \qquad\bigl(\mu\bigl(J_k^{\sigma_1
\sigma_2}\bigr)
\bigr)_{\sigma_1, \sigma_2
\in\{
0,1\}} = 2^{-k+1} e^{\sqrt{2} X(J_k) - \mathbb{E}X(J_k)^2 } \bigl( \mu
'\bigl(J_1^{\sigma_1 \sigma_2}\bigr) \bigr)_{\sigma_1, \sigma_2 \in\{
0,1\}},
\end{equation}
where $X(J_k) = W(\mathcal{C}(J_k))$ is a centered Gaussian random variable
with variance $\lambda(\mathcal{C}(J_k)) = (k-1) \log2$ and $\mu'$
is random
measure independent of $X(J_k)$ that has the same distribution as $\mu
$. But this implies that
\[
Z_k \stackrel{d} {=}2^{-2k+2} e^{2 \sqrt{2} X(J_k) - 2 \mathbb
{E}X(J_k)^2}
Z_1',
\]
where $Z_1' \stackrel{d}{=}Z_1$ is a random variable independent of
$Z_1$. Since
\[
2^{(-2k+2)h} \mathbb{E}e^{2 \sqrt{2} h X(J_k) - 2 h \mathbb
{E}X(J_k)^2} = 2^{(-2k+2)h} 2^{(4h^2-2h)(k-1)}
= 2^{4(h^2-h)(k-1)}
\]
and $\mathbb{E}Z_1^h$ is finite by Lemma \ref
{lemmamixedmoments-disjoint}, we have
\[
\mathbb{E}Z^h \leq\mathbb{E}Z_1^h \sum
_{k=1}^\infty2^{4(h^2-h)(k-1)} < \infty.
\]\upqed
\end{pf}

%
\begin{rem}\label{remfinitenessofc1}
While it can be seen from the proof of Proposition \ref
{propF-alpha-beta-asymptotics}, we emphasize that the finiteness of
$c_1$ follows from this lemma: simply use the elementary inequality
$\log(1+x)\leq\sqrt{x}$ for $x\geq0$ to bound $c_1$ by a term
proportional to $\mathbb{E}(\mu([0,\frac{1}{2}])^{{1}/{2}}\mu
([\frac
{1}{2},1])^{{1}/{2}})$.
\end{rem}

\begin{pf*}{Proof of Theorem \ref{teotail}}
We will show that for any $r > 1$ there exists a $\lambda_r$ such that
\[
\mathbb{P}(Y > \lambda) \leq\frac{c_1}{\lambda} r \qquad\mbox{for all }
\lambda
\geq\lambda_r.
\]
The verification of the lower bound is similar and is left to the reader.

Let $r > 1$ and fix $q > 1$ so that $\frac{q \log q}{q - 1} < \sqrt
{r}$. By Proposition \ref{propF-alpha-beta-asymptotics}, there exists
a $\lambda_r$ such that
\[
F_{1,q}(x) \leq c_1 \sqrt{r} \log q \qquad\mbox{for all }
x\geq\log\lambda_r,
\]
where we have defined $F_{1,q}(x) = \mathbb{E}(Y \mathbf{1}_{\{Y \in
{(e^x, q e^x]}\}
})$. We now have for $\lambda\ge\lambda_r$
\begin{eqnarray*}
\mathbb{P}( Y > \lambda) &=&  \sum_{k=0}^\infty
\mathbb{P} \bigl( Y \in\bigl( \lambda q^k, \lambda q^{k+1} \bigr]
\bigr)
\\
& \leq& \frac
{1}{\lambda
} \sum_{k=0}^\infty
q^{-k} F_{1,q}(k \log q + \log\lambda)
\\
& \leq& \frac{1}{\lambda} \sum_{k=0}^\infty
q^{-k} c_1 \sqrt{r} \log q = \frac{c_1}{\lambda} \sqrt{r}
\frac{q \log q}{q-1} \leq\frac{c_1}{\lambda} r,
\end{eqnarray*}
as was to be shown.
\end{pf*}

\section{Modulus of continuity}

\subsection{Outline of the proof}

In this section, we prove Theorem \ref{teomodulus}. Our plan of attack
is to follow the arguments carried out in \cite{bknsw12} in the case of
multiplicative cascades. However, the delicate dependence structure of
multiplicative chaos calls for nontrivial modifications. Let us briefly
sketch the main steps in the case of multiplicative cascades to see
what the main structure of the proof will be and what kind of
modifications we shall need.

The main part of the proof in the situation for cascades was showing
that if we write $(I_\sigma)_{\sigma\in\lbrace0,1\rbrace^n}$ for the
dyadic subintervals of $[0,1]$ of length $2^{-n}$ and $\mu$ for the
critical measure, then for any $\epsilon>0$ there exists a $C_\epsilon
>0$ such that for $\gamma\in(0,\frac{1}{2})$, $\mathbb{P}(\max_{\sigma
\in\lbrace0,1\rbrace^n}\mu(I_\sigma)\geq n^{-\gamma})\leq
C_\epsilon
n^{(1-\epsilon)(\gamma-({1}/{2}))}$. The corresponding result for the
modulus of continuity then follows from this through a Borel--Cantelli argument.

To\vspace*{1pt} get a hold of this estimate, one uses the scaling relation $(\mu
(I_\sigma))_\sigma\stackrel{d}{=}(e^{X_\sigma}Y^{(\sigma)})_\sigma$,
where $Y^{(\sigma)}$ are i.i.d. copies of $\mu([0,1])$ which are also
independent of the random variables $(X_\sigma)_\sigma$. By using the
scaling relation, conditioning on $(X_\sigma)$ and the tail estimate
$\mathbb{P}(Y^{(\sigma)}\geq\lambda)\approx C\lambda^{-1}$ (along with
some technical details to justify the approximations used)
\begin{eqnarray*}
\mathbb{P} \Bigl(\max_{\sigma\in\lbrace0,1\rbrace^n}\mu(I_\sigma)<
n^{-\gamma} \Bigr)&=& \mathbb{E} \biggl(\prod
_{\sigma\in\lbrace
0,1\rbrace^n}
\bigl(1-\mathbb{P}\bigl(Y^{(\sigma)}\geq n^{-\gamma}e^{-X_\sigma
}|(X_\sigma
)\bigr) \bigr)\biggr)
\\
&\approx&\mathbb{E} \biggl(\prod_{\sigma\in\lbrace0,1\rbrace
^n}
\bigl(1-Cn^{\gamma} e^{X_\sigma}\bigr) \biggr)
\\
&\approx&\mathbb{E}e^{-C n^{\gamma} \sum_{\sigma\in\lbrace
0,1\rbrace^n}
e^{X_\sigma}}.
\end{eqnarray*}

The last term we can write as $\phi_n(Cn^{\gamma-({1}/{2})})$, where
$\phi_n$ is the Laplace transform of the correctly normalized total
mass: $\phi_n(t)=\mathbb{E}e^{-t S_n}$, where $S_n=\sqrt{n}\sum_{\sigma
\in
\lbrace0,1\rbrace^n}e^{X_\sigma}$. One can then prove that for any
fixed $q\in(0,1)$, $\sup_n \mathbb{E}(S_n^q)<\infty$. Using this
and Markov's
inequality, one can show that for $q<1$, $1-\phi_n(t)\leq C_q t^q$ from
which one concludes that
%
%
\begin{eqnarray}
\mathbb{P} \Bigl(\max_{\sigma\in\lbrace0,1\rbrace^n}\mu(I_\sigma)\geq
n^{-\gamma} \Bigr)&\leq& 1-\phi_n\bigl(Cn^{\gamma-({1}/{2})}\bigr)
\nonumber\\[-8pt]\\[-8pt]
&\leq&
C_\epsilon n^{(1-\epsilon)(\gamma-({1}/{2}))}.\nonumber
\end{eqnarray}

While this sketch swept a lot of the technical details under the rug,
it still forms the back bone of the proof and one can see some of the
difficulties that will be present in the case of multiplicative chaos.
Let us consider some of the differences we can expect to be present in
the current context. First of all, if we manage to prove the same
estimate for the maximum of the measure of dyadic intervals, the
Borel--Cantelli argument will go through in a similar manner. The first
major difference is the scaling relation. For the exactly scale
invariant critical measure, one has a similar distributional relation:
$(\mu_{\sqrt{2}}(I_\sigma))_\sigma\stackrel{d}{=}(e^{X_\sigma}\mu
^{(\sigma)}([0,1]))_\sigma$, but the difference is that we have
nontrivial correlations---$\mu^{(\sigma)}$ are not independent from
each other and they may depend on some of the $X_\sigma$ as well. To
remedy this, we consider instead of $\mu_{\sqrt{2}}$ another random
measure which is absolutely continuous with respect to $\mu_{\sqrt{2}}$ which possesses nice scaling
properties, nice decorrelation properties as well as a nicely behaving
Radon--Nikodym derivative with respect to~$\mu_{\sqrt{2}}$. Moreover,
one gets similar asymptotic behavior for the tail of the measure of the
unit interval for this measure as well.

The next step of the proof is to use scaling, independence and tail
behavior to obtain a similar estimate in terms of a Laplace transform
and some errors due to approximations. This step of the proof requires
a fair amount of technical details which are even more involved than in
the multiplicative cascade setup, but philosophically similar. Finally,
we are left with estimating moments of the correctly normalized
approximation to the critical measure. This can be done by using
Gaussian comparison inequalities and the result from multiplicative cascades.

\subsection{Tools for the proof} Let us now collect some of the tools
we shall need for the proof. First of all, we shall consider
modifications of the field $X$ and the measure $\mu_{\sqrt{2}}$ for
which we still have a similar result for the tail.

%
\begin{lemma}\label{letailmodification}
Assume\vspace*{1pt} that we can write $\mu_{\sqrt{2}}(\dd x)=e^{Z(x)}\nu(\dd x)$, for some
random measure $\nu(\dd x)$ and random Gaussian field $Z$ which is
independent of $\nu$ and $\min_{x\in[0,1]}Z(x)>0$ with positive
probability, then there exists a constant $C$ such that $\mathbb
{P}(\nu
([0,\alpha])>\lambda)\leq C \alpha\lambda^{-1}$.
\end{lemma}

\begin{pf}
Plugging in the definitions,
\begin{eqnarray*}
\mathbb{P} \bigl(\mu_{\sqrt{2}}\bigl([0,\alpha]\bigr)>\lambda\bigr)&=&
\mathbb{P} \biggl(\int_0^\alpha e^{Z(x)}
\nu(\dd x)>\lambda\biggr)
\\
&\geq&\mathbb{P} \bigl(e^{\min_{x\in[0,1]}Z(x)}\nu\bigl([0,\alpha]\bigr
)>\lambda\bigr)
\\
&\geq&\mathbb{P} \bigl(e^{\min_{x\in[0,1]}Z(x)}>1,\nu\bigl([0,\alpha
]\bigr)>\lambda
\bigr)
\\
&=&\mathbb{P} \Bigl(\min_{x\in[0,1]}Z(x)>0 \Bigr)\mathbb{P} \bigl(
\nu\bigl([0,\alpha]\bigr)>\lambda\bigr).
\end{eqnarray*}
On the other hand, by scaling
%
%
\begin{equation}
\mathbb{P} \bigl(\mu_{\sqrt{2}}\bigl([0,\alpha]\bigr)>\lambda\bigr)=
\mathbb{P} \bigl(\alpha e^{X_\alpha-({1}/{2})\mathbb{E}(X_\alpha
^2)}\mu_{\sqrt
{2}}\bigl([0,1]\bigr)>
\lambda\bigr),
\end{equation}
where $X_\alpha$ is a centered Gaussian independent of $\mu
_{\sqrt{2}}([0,1])$. Conditioning on $X_\alpha$ and using the tail
estimate for $\mu_{\sqrt{2}}([0,1])$
%
%
\begin{equation}
\mathbb{P}\bigl(\alpha e^{X_\alpha-({1}/{2})\mathbb{E}(X_\alpha
^2)}\mu_{\sqrt
{2}}\bigl([0,1]\bigr)>
\lambda\bigr)\leq C \alpha\lambda^{-1}.
\end{equation}

Collecting everything gives the desired result.
\end{pf}

%
\begin{rem}\label{remtailmodification}
While the class of measures $\nu$ covered by this result is rather
limited (due to the fact that the result was easy to prove and
sufficient for our needs concerning the modulus of continuity), we
believe that such a result for the tail should hold quite generally for
critical Gaussian multiplicative chaos measures.
\end{rem}

We next note that the regular variation with exponent $-1$ of the tail
is robust under linear combinations of copies of random variables:

%
\begin{lemma}\label{lemmasum-inequality}
Let\vspace*{1pt} $X\geq0$ satisfy $\mathbb{P}(X>\lambda)\leq\frac
{A}{\lambda}$ for $\lambda> 0$.

Let $X_j$, $j\in\lbrace1,\ldots,N\rbrace$ be ({possibly
dependent}) random variables with the
same distribution as $X$ and
let $a_j\geq0$ for $j\in\lbrace1,\ldots,N\rbrace$. Then
\[
\mathbb{P} \Biggl(\sum_{j=1}^N
a_j X_j>\lambda\Biggr)\leq\frac{C\cdot
A\log
(N+1) (\sum_{j=1}^N a_j )}{\lambda}\qquad\mbox{for all } \lambda>0,
\]
with a universal ({in particular, independent of $A$})
constant $C<\infty$.
\end{lemma}

\begin{pf}
We may assume that $\sum_{j=1}^Na_j=1$ since
the statement scales in the right way. Fix $t\in(0,1)$ and observe
first that for all positive $y_1,\ldots,y_N$
one has the subadditivity inequality
\[
\Biggl(\sum_{j=1}^Na_j
y_j \Biggr)^t\leq\sum_{j=1}^Na_j^ty_j^t.
\]
Fix $\lambda>0$.
The above holds if we set $y_j=(x_j-\lambda)_+$, where we denote the positive
part by $y_+:=\max(0,y)$ and let, for now, the numbers $(x_j)_{1\leq j
\leq N}$ be arbitrary reals.
We obtain using $\sum_{j=1}^Na_j=1$ (and Jensen) that
\[
\Biggl(\sum_{j=1}^Na_j
x_j-\lambda\Biggr)_+^t \leq\Biggl(\sum
_{j=1}^Na_j (x_j-
\lambda)_+ \Biggr)^t \leq\sum_{j=1}^Na_j^t(x_j-
\lambda)_+^t,
\]
or, in other words,
\[
\phi\Biggl(\sum_{j=1}^Na_j
x_j\Biggr)\leq\sum_{j=1}^Na_j^t
\phi(x_j),
\]
where $\phi(x):=(x-\lambda)_+^t$. Especially, we have
%
%
\begin{equation}
\label{0} \mathbb{E}\phi\Biggl(\sum_{j=1}^N
a_j X_j\Biggr) \leq\mathbb{E}\phi(X)\sum
_{j=1}^N a_j^t.
\end{equation}
The right-hand side can be estimated as follows:
%
%
\begin{eqnarray}
\label{1} \mathbb{E}\phi(X) &=& \int_{0}^\infty
\phi'(u)\mathbb{P}(X>u) \,\dd u\nonumber
\\
& \leq& A\int_{\lambda}^\infty
t(u-\lambda)^{t-1}u^{-1} \,\dd u
\nonumber
\\
&=&At\int_{0}^\lambda\frac{y^{t-1}\,\dd y}{y+\lambda}+ At\int
_{\lambda}^\infty\frac{y^{t-1}\,\dd y}{y+\lambda}
\\
&\leq& At\lambda^{-1}\int_{0}^\lambda
y^{t-1}\,\dd y+ At\int_{\lambda}^\infty
y^{t-2}\,\dd y
\nonumber
\\
&=&A(1-t)^{-1}\lambda^{t-1}.\nonumber
\end{eqnarray}
From Markov's inequality and (\ref{0}), we thus obtain
\[
\phi(2\lambda) \cdot\mathbb{P}\Biggl(\sum_{j=1}^N
a_j X_j>2\lambda\Biggr) \leq\mathbb{E}\phi\Biggl(
\sum_{j=1}^N a_j
X_j\Biggr) \leq\Biggl(\sum_{j=1}^Na_j^t
\Biggr)\mathbb{E}\phi(X),
\]
and by combining with (\ref{1})
%
%
\begin{equation}
\label{23} \mathbb{P}\Biggl(\sum_{j=1}^N
a_j X_j>2\lambda\Biggr)\leq\frac{A}{\lambda
} \Biggl(
\frac{1}{1-t}\sum_{j=1}^Na_j^t
\Biggr).
\end{equation}
Finally, choosing $t=t_0:=1-1/\log N$ (for $N\geq3$) we get
\[
\Biggl(\frac{1}{1-t_0}\sum_{j=1}^Na_j^{t_0}
\Biggr) \leq\biggl(\frac
{N^{1-t_0}}{1-t_0} \biggr) \Biggl(\sum
_{j=1}^Na_j \Biggr)^{t_0} =
\frac
{N^{1-t_0}}{1-t_0} = e \log N,
\]
and then (\ref{23}) yields the stated result.
\end{pf}

%
\begin{rem}
The above result is essentially optimal: choose $\Omega= [0,1)$, that
is, the one-dimensional torus with the Lebesgue measure.
Let
\[
X_0(\omega)=\frac{N}{k} \qquad\mbox{for }\omega\in
\bigl[(k-1)/N,k/N\bigr),  k=1,2,\ldots,N.
\]
Then
$\mathbb{P}(X>\lambda)<1/\lambda$. Define
the random variables $X_j$, $j=1,\ldots,N$ with the formula
$X_j(\omega)=X_0(\omega+(j-1)/N)$, which is
well defined since we are now in the torus. Then each $X_j$ has the
same tail as $X_0$. However, the average $X:=(1/N)\sum_{j=1}^NX_j$
is the constant\vspace*{1pt} variable: $X(\omega)=\sum_{j=1}^Nj^{-1}\geq\log N$ for
all $\omega\in\Omega$. We thus have
$\mathbb{P}(X\geq\log N)=1$.
\end{rem}

For comparing the present setting with that of multiplicative cascades, we shall make use
of Kahane's convexity inequalities \cite{K2}.

%
\begin{prop}\label{propkahane-inequality}
Let $G\dvtx [0,\infty)\rightarrow\mathbb{R}$ be a concave function such that
$|G(x)|\leq C(1+x^\alpha)$ for some positive constants $C$ and $\alpha
$. Let $A \subset\mathbb{R}^d$ be a Borel set, $\rho$ be a Radon
measure on
$A$ and $(X_r)_{r\in A}$ and $(Y_r)_{r\in A}$ be two continuous and
centered Gaussian processes on $A$ such that the covariance kernels
satisfy $k_X(u,v)\leq k_Y(u,v)$ for all $u,v\in A$. Then
\[
\mathbb{E} G \biggl(\int_A e^{X_r-({1}/{2})\mathbb
{E}(X_r^2)}\rho(\d r)
\biggr) \geq\mathbb{E} G \biggl(\int_A e^{Y_r-
({1}/{2})\mathbb
{E}(Y_r^2)}
\rho(\d r) \biggr).
\]
\end{prop}

To apply this inequality, we construct a Gaussian field on $[0,1]$ for
which the moments of the corresponding measure can be calculated and
for which we have a covariance structure that allows comparing with
more correlated situations (such a comparison is also used in \cite
{drsv12-1} to prove that the limit of the total mass martingale
associated to nonrenormalized critical chaos measures vanishes almost surely).

The Gaussian field we shall employ is essentially a Gaussian branching
random walk. Let us associate to the collection $\{I_\sigma\}$ of
dyadic subintervals of $[0,1]$ an i.i.d. collection of standard
Gaussian random variables $\{V_\sigma\}$. Let us write $\Sigma
_k=\lbrace0,1\rbrace^k$ and define the field
%
%
\begin{equation}
\label{un} U_n(x)=\sum_{k=1}^n
\sum_{\sigma\in\Sigma_k\dvtx  x\in I_\sigma} V_\sigma.
\end{equation}
The covariance of $U_n$ is given by
\begin{eqnarray*}
\mathbb{E}\bigl(U_n(x)U_n(y)\bigr)&=&\sum
_{k,k'=1}^n\sum_{\sigma\in\Sigma
_k\sigma
'\in\Sigma_{k'}\dvtx  x\in I_\sigma,y\in I_{\sigma'}}
\mathbb{E}(V_\sigma V_{\sigma'})
\\
&=&\sum_{k,k'=1}^n \sum
_{\sigma\in\Sigma_k,\sigma'\in\Sigma
_{k'}\dvtx  x\in
I_\sigma, y\in I_{\sigma'}}\mathbf{1}\bigl(\sigma=\sigma'\bigr)
\\
&=&\sum_{k=1}^n \sum
_{\sigma\in\Sigma_k\dvtx  x,y\in I_\sigma} 1.
\end{eqnarray*}

For comparison with other fields, we note that to have a $\sigma\in
\Sigma_k$ such that \mbox{$x,y\in I_\sigma$,} we must have $|x-y|\leq2^{-k}$
and we see that
\begin{eqnarray*}
\mathbb{E}\bigl(U_n(x)U_n(y)\bigr)&\leq&\sum
_{k=1}^{({-\log|x-y|}/{\log2})
\wedge n} 1
\\
&=&\frac{-\log|x-y|}{\log2} \wedge n.
\end{eqnarray*}

Our last technical lemma is a version of the
Borell--Tsirelson--Ibragimov--Sudakov inequality \cite{adta07}, Theorem
2.1.1. For our purposes, we need a version which relates the
tail probability of the supremum of a Gaussian process on an interval
both to the size of the interval and to the modulus of continuity of
the covariance of the process in a quantitative manner.

%
\begin{lemma}\label{lemmaborell-tsirelson}
Let\vspace*{1pt} $I \subset\mathbb{R}$ be a bounded interval and $L > 0$. Let $(Y(x))_{x
\in I}$ be an arbitrary centered Gaussian process on $I$ such that
$\mathbb{E}
|Y(x)-Y(y)|^2 \leq L|x-y|$ for all $x,y \in I$, and further suppose
there is some (deterministic) $x_0 \in I$ for which $Y(x_0) = 0$ almost
surely. Then, for any $\varepsilon> 0$, there exists an absolute
constant $c_\varepsilon> 0$ (i.e., the choice of $c_\varepsilon$
depends only on $\varepsilon$) such that for all $s > 0$
%
%
\begin{equation}
\label{eqborell-tsirelson} \mathbb{P} \Bigl( \sup_{x \in I} Y(x) > s
\Bigr)
\leq c_\varepsilon e^{-
{s^2}/{((2+\varepsilon) |I| L)}}.
\end{equation}
\end{lemma}

\begin{pf}
By considering the scaled process $\frac{1}{\sqrt{|I| L}} Y(|I|\cdot)$
instead of $Y(\cdot)$ we may without loss of generality reduce to the
case $|I| = L = 1$. Since $\mathbb{E}Y(x_0)^2 = 0$, this normalization also
implies that $\sigma_Y^2:= \sup_{x \in I} \mathbb{E}Y(x)^2 \leq1$. The
Borell--TIS inequality then states that for $s > 0$ we have
%
%
\begin{equation}
\label{eqborell-tsirelson-proof} \mathbb{P} \Bigl( \sup_{x \in I} Y(x)
- \mathbb{E}
\sup_{x \in I} Y(x) > s \Bigr) \leq e^{-{s^2}/{(2 \sigma_Y^2)}} \leq
e^{-{s^2}/{2}}.
\end{equation}
Then consider the Gaussian process $X(x) = B_x-B_{x_0}$, where
$(B_x)_{x \in I}$ is a one-dimensional Brownian motion. Clearly,
$(X(x))_{x \in I}$ satisfies the assumptions of the lemma, and moreover,
\[
\mathbb{E}\bigl|Y(x)-Y(y)\bigr|^2 \leq|x-y| = \mathbb{E}\bigl|X(x)-X(y)\bigr|^2
\]
for all $x,y \in I$. By the Sudakov--Fernique inequality (\cite{adta07}, Theorem
2.2.3), we then have
\[
\mathbb{E}\sup_{x \in I} Y(x) \leq\mathbb{E}\sup
_{x \in I} X(x) \leq M < \infty
\]
for some absolute constant $M > 0$. In (\ref
{eqborell-tsirelson-proof}), for $s > M$ this implies
\[
\mathbb{P} \Bigl( \sup_{x \in I} Y(x) > s \Bigr) \leq
e^{-{(s-M)^2}/{2}}.
\]
Since the choice of $M$ does not depend on the parameters of the
process $(Y(x))_{x \in I}$, it clear that for any $\varepsilon> 0$
there exists an absolute constant $c_\varepsilon> 0$ for which (\ref
{eqborell-tsirelson}) holds.
\end{pf}

%
\begin{rem}
The statement of the lemma generalizes to processes on bounded domains
$U \subset\mathbb{R}^d$ for $d \geq2$ simply by replacing the length
$|I|$ of
the interval $I$ by the diameter $\operatorname{diam}(U)$ of $U$. The only
difference in the proof is that instead of one-dimensional Brownian
motion one compares the arbitrary process to L\'{e}vy's Brownian motion
on $\mathbb{R}^d$, that is, the Gaussian process $(X(x))_{x \in
\mathbb{R}^d}$ with $\mathbb{E}
X(x) X(y) = \frac{1}{2}(|x|+|y|-|x-y|)$; we refer to \cite{K4} for a
proof that this function is indeed a covariance kernel.
\end{rem}

We are ready to proceed to the main proof.

\subsection{Main results for the modulus of continuity}

Let $((X_t(x))_{x\in\mathbb{R}})_{t\geq0}$ be the exactly scale invariant
Gaussian field on $\mathbb{R}$ as before and define the Gaussian field
$((Y_t(x))_{x\in\mathbb{R}})_{t\geq0}$ by setting
\[
Y_t(x) = W\bigl(\mathcal{C}_t(x) \setminus
\mathcal{C}_0(x)\bigr) = X_t(x) - X_0(x)
\qquad\mbox{for } x \in\mathbb{R}, t \geq0.
\]
In the proof of Theorem \ref{teomodulus}, it is convenient to use the
characterization (\ref{eqrenormalized-critical-measure}) of critical
lognormal multiplicative chaos. To keep the notation simpler, we
normalize the construction by the deterministic constant $c > 0$ in
(\ref{eqrenormalized-critical-measure}). Explicitly, we consider the
critical measures associated to the fields $X$ and $Y$ and denote
\[
\mu_{\sqrt{2}}(\d x)=\lim_{t\rightarrow\infty}\sqrt{t}e^{\sqrt
{2}X_t(x)-(t+1)}
\,\d x
\]
and
%
\[
\nu_{\sqrt{2}}(\d x)=\lim
_{t\rightarrow\infty}\sqrt{t}e^{\sqrt
{2}Y_t(x)-t}\,\d x,
\]
where the limits exist in probability in the weak sense. By
construction, it is clear that almost surely, the Radon--Nikodym
derivative $\frac{\d\mu_{\sqrt{2}}}{\d\nu_{\sqrt{2}}}(x) =
e^{\sqrt
{2} X_0(x) - 1}$ is almost surely positive and uniformly bounded away
from $0$ and $\infty$ for all $x \in[0,1]$ (in particular, the
assumptions of Lemma \ref{letailmodification} are met), so for the
purpose of our result on the modulus of continuity the difference
between these two measures is insignificant. The measure $\mu_{\sqrt
{2}}$ is exactly scale invariant as before, but in this section we make
more use of the measure $\nu_{\sqrt{2}}$ which satisfies the \emph
{$\star$-scaling relation}: for every $\epsilon\in(0,1]$ we have
%
%
\begin{equation}
\label{eqstar-scaling} \bigl(\nu_{\sqrt{2}}(A)\bigr)_{A\in\mathcal
{B}([0,1])}\stackrel{d} {=}
\biggl(\epsilon\int_A e^{\sqrt{2}Y_{-\log\epsilon}+\log\epsilon}\nu
_{\sqrt
{2}}^\epsilon(\d x) \biggr)_{A\in\mathcal{B}([0,1])},
\end{equation}
where $\nu_{\sqrt{2}}^\epsilon$ is independent of $Y_{-\log\epsilon}$
and $(\nu_{\sqrt{2}}^\epsilon(A))_{A}\stackrel{d}{=}(\nu_{\sqrt
{2}}(\epsilon^{-1}A))_{A}$. The proof of this scaling relation is
recalled in the \hyperref[app]{Appendix}. We also stress that $\nu_{\sqrt{2}}$
satisfies the conditions of Lemma \ref{letailmodification}.

The next lemma contains the key technical estimates that lead to the
proof of Theorem \ref{teomodulus}.

%
\begin{lemma}\label{lemmamodulus-probability}
Let us index by $\sigma\in\Sigma_n=\lbrace0,1\rbrace^n$ the dyadic
subintervals $I_\sigma$ of $[0,1]$ of length $2^{-n}$. Moreover, write
$\Sigma_n^{(e)}$ for the family of even dyadic intervals of length
$2^{-n}$ (i.e., intervals of the form $[(2j)2^{-n},(2j+1)2^{-n})$). Then
for $\gamma\in(0,\frac{1}{2})$ and $\epsilon\in(0,1)$ there exists a
constant $C=C(\epsilon)$ such that
%
%
\begin{equation}
\mathbb{P}\Bigl(\max_{\sigma\in\Sigma_n^{(e)}}\nu_{\sqrt{2}}(I_\sigma
)\geq n^{-\gamma}\Bigr)\leq C n^{(1-\epsilon)(\gamma-({1}/{2}))}.
\end{equation}
The same holds if we replace $\Sigma_n^{(e)}$ with $\Sigma_n^{(o)}$,
the corresponding collection of odd dyadic intervals.
\end{lemma}

\begin{pf}
The proof is rather lengthy so we shall split it into steps that
somewhat parallel the cascade proof.
\begin{longlist}[\textit{Step} 2:]
\item[\textit{Step} 1:] \textit{Using scaling and independence}.

We begin by noting that by specializing the $\star$-scaling relation to
dyadics, we get
%
%
\begin{equation}
\bigl(\nu_{\sqrt{2}}(I_\sigma)\bigr)_{\sigma\in\Sigma_n}\stackrel{d}
{=} \biggl( 2^{-n}\int_{I_\sigma} e^{\sqrt{2}Y_{n\log2}(x)-n\log2}
\nu_{\sqrt
{2}}^{(n)}(\d x) \biggr)_{\sigma\in\Sigma_n},
\end{equation}
where $\nu_{\sqrt{2}}^{(n)}$ is independent of $Y_{n\log2}$ and
$(\nu
_{\sqrt{2}}^{(n)}(A))_A\stackrel{d}{=}(\nu_{\sqrt{2}}(2^n A))_A$. Since
$Y_t(x)$ and $Y_t(y)$ are independent when $|x-y|\geq1$, $\nu_{\sqrt
{2}}(A)$ is independent\break of $\nu_{\sqrt{2}}(B)$ when $d(A,B)\geq1$.
Thus, the scaling property implies that\break $(\nu_{\sqrt{2}}^{(n)}\lfloor
I_\sigma)_{\sigma\in\Sigma_n^{(e)}}$ is a family of independent random
measures (and similarly for the odd intervals)---here $\nu_{\sqrt
{2}}^{(n)}\lfloor I_\sigma$ denotes the restriction of $\nu_{\sqrt
{2}}^{(n)}$ to $I_\sigma$.

Let us write
%
%
\begin{equation}
\label{eqWn} W_{n,\sigma}=2^{-n}\int_{I_\sigma}
e^{\sqrt{2}Y_{n\log2}(x)-n\log
2}\nu_{\sqrt{2}}^{(n)}(\d x).
\end{equation}
Using the independence noted above, we see that
\begin{eqnarray}\label{ineq1}
\nonumber
\mathbb{P} \Bigl(\max_{\sigma\in\Sigma_n^{(e)}} W_{n,\sigma
}<n^{-\gamma}
\Bigr) &=& \mathbb{E} \prod_{\sigma\in\Sigma_n^{(e)}} \mathbb{P}
\bigl(W_{n,\sigma}<n^{-\gamma}|Y_{n\log2}\bigr)
\nonumber\\[-8pt]\\[-8pt]
& \geq&\mathbb{E} \prod_{\sigma\in\Sigma
_n}\bigl(1-
\mathbb{P}\bigl(W_{n,\sigma}\geq n^{-\gamma}|Y_{n\log2}\bigr)
\bigr).\nonumber
\end{eqnarray}

\item[\textit{Step} 2:] \textit{Getting to the Laplace transform}.

To estimate $\mathbb{P}(W_{n,\sigma}\geq n^{-\gamma}|Y_{n\log2})$, we
will approximate the integral (\ref{eqWn}) by a Riemann sum and then
make use of Lemma \ref{letailmodification}. 
For brevity, we will denote $f_{\sqrt{2}}^{(n)}(\cdot):=e^{\sqrt
{2}Y_{n\log2}(\cdot)-n\log2}$. Fix $\sigma\in\Sigma_n$ for the
moment, let $k \in\mathbb{N}_+$ and divide $I_\sigma$ into $2^k$ subintervals
$(I_{\sigma,j})_{j=1}^{2^k}$ of equal length. Denote the midpoint of
$I_{\sigma,j}$ by $x_{\sigma,j}$. Let $s>0$ and define the event
$\mathcal{D}_s= \lbrace\sup_{x \in I_{\sigma,j}} |Y_{n \log
2}(x) -
Y_{n \log2}(x_{\sigma,j})| \leq s$ for all $j=1,2,\ldots,2^k
\rbrace
$. We then have on $\mathcal{D}_s$
\[
2^n W_{n,\sigma} = \int_{I_\sigma}
f_{\sqrt{2}}^{(n)}(x) \nu_{\sqrt
{2}}^{(n)}(\d x)
\leq e^{2\sqrt{2} s} \sum_{j=1}^{2^k}
\Xint-_{I_{\sigma,j}} f_{\sqrt{2}}^{(n)}(x) \,\d x\,
\nu_{\sqrt{2}}^{(n)}(I_{\sigma,j}),
\]
where\vspace*{-2pt} $\Xint-_A f(x) \,\d x:= \frac{1}{|A|}\int_A f(x) \,\d x$ is the
integral average. Let $\mathcal{F}_n=\sigma(\{Y_t(x)\dvtx\break   x\in[0,1],
t\leq n\log2\})$. Since $\nu_{\sqrt{2}}^{(n)}(I_{\sigma,j})
\stackrel{d}{=}\nu
_{\sqrt{2}}^{(n)}(I_{\sigma,i})$ for $j \neq i$ and the function
$f_{\sqrt{2}}^{(n)}$ is independent of the measure $\nu_{\sqrt
{2}}^{(n)}$, Lemmas \ref{lemmasum-inequality}~and~\ref
{letailmodification} imply that on~$\mathcal{D}_s$
\begin{eqnarray*}
&& \mathbb{P} \biggl( \int_{I_\sigma} f_{\sqrt{2}}^{(n)}(x)
\nu_{\sqrt
{2}}^{(n)}(\d x) > \lambda\Big\vert\mathcal{F}_n
\biggr)
\\
&&\qquad  \leq
\mathbb{P} \Biggl( e^{2\sqrt{2} s} \sum_{j=1}^{2^k}
\Xint-_{I_{\sigma,j}} f_{\sqrt{2}}^{(n)}(x) \,\d x\,
\nu_{\sqrt
{2}}^{(n)}(I_{\sigma,j}) > \lambda\Big\vert
\mathcal{F}_n \Biggr)
\\
&&\qquad \leq Ck 2^{-k} \biggl( \frac{e^{2\sqrt{2} s} \sum_{j=1}^{2^k} \Xint
-_{I_{\sigma,j}} f_{\sqrt{2}}^{(n)}(x) \,\d
x}{\lambda} \biggr)
\end{eqnarray*}
for some constant $C>0$. Setting $\lambda=n^{-\gamma} 2^n$ and
combining this inequality with (\ref{ineq1}) and the inequality
$e^{-2x} \leq1-x$ valid for $x\in[0,1/2]$, we get
\begin{eqnarray*}
&& \mathbb{E}\mathbb{P} \Bigl( \max_{\sigma\in\Sigma_n^{(e)}} \nu
_{\sqrt{2}}(I_\sigma) > n^{-\gamma} \big\vert
\mathcal{F}_n \Bigr)
\\
&&\qquad \leq1 - \mathbb{E}\prod_{\sigma\in\Sigma
_n}
\bigl(1-\mathbb{P}\bigl(W_{n,\sigma}\geq n^{-\gamma}|Y_{n\log2}
\bigr)\bigr)
\\
&&\qquad \leq1 - \mathbb{E}\exp\biggl( - 2 C k 2^{-k}
e^{2\sqrt{2} s} \sum_{\sigma\in\Sigma_n} \biggl(
\frac{\sum_{j=1}^{2^k} \Xint-
_{I_{\sigma,j}} f_{\sqrt{2}}^{(n)}(x) \,\d x}{2^{n} n^{-\gamma}} \biggr)
\biggr)
\\
&&\quad\qquad{}+ 1 - \mathbb{P}( \mathcal{A}_{n,k,s} )
\\
&&\qquad = 1-\mathbb{E}\exp\biggl(-2Ck e^{2\sqrt{2}s}
n^\gamma\int_0^1 f_{\sqrt{2}}^{(n)}(x)
\,\d x \biggr)+1-\mathbb{P}(\mathcal{A}_{n,k,s}),
\end{eqnarray*}
where $\mathcal{A}_{n,k,s}$ is the event
\begin{eqnarray*}
\mathcal{A}_{n,k,s} &=& \biggl\{\max_{\sigma\in\Sigma_n} Ck
2^{-k} \frac
{e^{2\sqrt{2} s} \sum_{j=1}^{2^k} \Xint-_{I_{\sigma,j}} f_{\sqrt
{2}}^{(n)}(x) \,\d x}{n^{-\gamma} 2^n}<\frac{1}{2} \biggr\}
\\
&&{} \cap\Bigl\{ \sup_{x \in I_{\sigma,j}} \bigl|Y_{n \log2}(x) -
Y_{n
\log2}(x_{\sigma,j})\bigr| \leq s
\\
&&\hspace*{31pt} \forall j \in\bigl\{0,1,
\ldots,2^k-1\bigr\}\ \forall\sigma\in\Sigma_n \Bigr\}.
\end{eqnarray*}
%
Denoting
\[
S_{n} = n^{{1}/{2} } \int_0^1
e^{\sqrt{2} Y_{n \log2}(x) - n
\log
2} \,\d x,
\]
we finally get
\begin{eqnarray}\label{eqmaximum-prob-estimate}
\mathbb{P} \Bigl( \max_{\sigma\in\Sigma_n^{(e)}} \mu_{\sqrt
{2}}(I_\sigma)
> n^{-\gamma} \Bigr) & \leq&1 - \mathbb{E}\exp\bigl( - 2Ce^{2\sqrt{2}s}k
n^{ (\gamma-({1}/{2}) )} S_{n} \bigr)
\nonumber
\nonumber\\[-8pt]\\[-8pt]
&&{} + 1 - \mathbb{P}(\mathcal{A}_{n,k,s}). \nonumber
\end{eqnarray}

\item[\textit{Step} 3:] \textit{Controlling the error}.

We then estimate the terms in the inequality above. Denote
\begin{eqnarray*}
\mathcal{B}_n &=& \biggl\{\max_{\sigma\in\Sigma_n} Ck
2^{-k} \frac
{e^{2\sqrt{2} s} \sum_{j=1}^{2^k} \Xint-_{I_{\sigma,j}} f_{\sqrt
{2}}^{(n)}(x) \,\d x}{n^{-\gamma} 2^n}<\frac{1}{2} \biggr\}
\\
&=& \biggl\lbrace\max_{\sigma\in\Sigma_n}\int_{I_\sigma}
e^{\sqrt{2}
Y_{n\log2}(x)-n\log2}\,\d x< n^{-\gamma}\bigl(2Ck e^{2\sqrt
{2}s}
\bigr)^{-1} \biggr\rbrace
\end{eqnarray*}
and
\[
\mathcal{B}_{n,k,s}'= \Bigl\{ \sup_{x \in I_{\sigma,j}}
\bigl|Y_{n \log
2}(x) - Y_{n \log2}(x_{\sigma,j})\bigr| \leq s \ \forall j
\in\bigl\{0,1,\ldots,2^k-1\bigr\}\ \forall\sigma\in
\Sigma_n \Bigr\}
\]
so that
\[
\mathcal{A}_{n,k,s} = \mathcal{B}_n \cap
\mathcal{B}_{n,k,s}' \quad\mbox{and} \quad1 - \mathbb{P}(
\mathcal{A}_{n,k,s}) \leq\bigl(1 - \mathbb{P}(\mathcal
{B}_n) \bigr) + \bigl(1 - \mathbb{P}\bigl(\mathcal{B}_{n,k,s}'
\bigr) \bigr).
\]
We first estimate the probability of $\mathcal{B}_{n,k,s}'$ not occurring.
For all $\sigma$ and $j$, the length of $I_{\sigma,j}$ is $2^{-n-k}$
and $\mathbb{E}|Y_{n \log2}(x)-Y_{n \log2}(y)|^2 \leq2^{n+1} |x-y|$,
so by Lemma \ref{lemmaborell-tsirelson} we have, for any $\sigma
\in
\Sigma_n$ and $j=1,\ldots,2^k$,
\[
\mathbb{P} \Bigl( \sup_{x \in I_{\sigma,j}} \bigl|Y_{n \log2}(x) -
Y_{n
\log
2}(x_{\sigma,j})\bigr| > s \Bigr) \leq c e^{-2^{k-3} s^2},
\]
where $c > 0$ is an absolute constant. It follows that
\[
1 - \mathbb{P}\bigl(\mathcal{B}_{n,k,s}'\bigr) \leq c
2^{n+k} e^{-2^{k-3} s^2}.
\]
For the choice $s_n \sim\sqrt{\epsilon\log n}$ and $k_n \sim\alpha
\log n$, the right-hand side of this estimate is asymptotically
equivalent to
\[
n^{\alpha\log2} e^{n \log2 - ({\epsilon}/{8}) n^{\alpha\log2}
\log n},
\]
from\vspace*{2pt} which we see that in order to have $\sum_{n=1}^\infty(1-\mathbb{P}
(\mathcal{B}_{n,k_n,s_n}')) < \infty$ we may take $\epsilon> 0$ arbitrarily
small, but must restrict to $\alpha\geq1/\log2$. Taking $\alpha=
1/\log2$, in~(\ref{eqmaximum-prob-estimate}) these choices give
%
%
\begin{eqnarray}\label{eqPmax}
&& \mathbb{P} \Bigl( \max_{\sigma\in\Sigma_n^{(e)}}
\nu_{\sqrt
{2}}(I_\sigma) > n^{-\gamma} \Bigr)\nonumber
\\
&&\qquad  \leq  1 -
\mathbb{E}\exp\biggl( -2 C e^{2\sqrt
{2}\sqrt
{\epsilon\log n}} \frac{\log n}{\log2} n^{ (\gamma-({1}/{2}) )}
S_{n} \biggr)
\\
&&\quad\qquad{} + c' n^{-c'' \log n} + \bigl(1 - \mathbb{P}(
\mathcal{B}_n) \bigr)\nonumber
\end{eqnarray}
for some constants $c', c'' > 0$ depending on $\epsilon$.

To estimate the probability of $\mathcal{B}_n$, we note that
%
%
\begin{equation}
\bigl\lbrace S_n<n^{({1}/{2})-\gamma}\bigl(2Ck_n
e^{2\sqrt
{2}s_n}\bigr)^{-1}\bigr\rbrace\subset\mathcal{B}_n.
\end{equation}
By Chebyshev's inequality, we then see that for any $q<1$
\begin{eqnarray*}
1-\mathbb{P}(\mathcal{B}_{n}) & \leq&\mathbb{P}\bigl(S_{n}>
\bigl(2Ck_n e^{2\sqrt
{2}s_n}\bigr)^{-1}n^{(({1}/{2})-\gamma)}
\bigr)
\\
& \leq&\bigl(2Ck_n e^{2\sqrt{2}s_n}\bigr)^{q}
\frac{\mathbb
{E}(S_{n}^q)}{n^{(({1}/{2})-\gamma)q}}.
\end{eqnarray*}

\item[\textit{Step} 4:] \textit{Comparison with cascades}.

If we knew that $\mathbb{E}(S_{n}^q)$ were uniformly bounded in $n$ for
some values of $q$, we would have a quantitative estimate for the speed
at which $\mathbb{P}(\mathcal{B}_{n})$ tends to one. For this, we
employ Kahane's convexity inequalities, that is, Proposition \ref
{propkahane-inequality}, and comparison with the branching random walk
$U_n$ defined in (\ref{un}).

Note that
\begin{eqnarray*}
\mathbb{E}\bigl(U_n(x)U_n(y)\bigr)&\leq&
\frac{-\log|x-y|}{\log2} \wedge n
\\
&\leq&\frac{1}{\log2}\mathbb{E}\bigl(Y_{n\log2}(x)Y_{n\log2}(y)
\bigr)+C,
\end{eqnarray*}
for some large enough constant $C$, since the covariance of the field
$Y_{n\log2}$ is given by
\begin{eqnarray*}
&& \mathbb{E}\bigl(Y_{n\log2}(x)Y_{n\log2}(y)\bigr)
\\
&&\qquad = \cases{
-\log|x-y|+|x-y|-1, &\quad$2^{-n}\leq|x-y|\leq1$,
\vspace*{3pt}\cr
n
\log2+|x-y|-2^n|x-y|, &\quad$ |x-y|\leq2^{-n}$.}
\end{eqnarray*}

Let us thus consider a standard Gaussian variable $Z$ independent of
$Y_{n\log2}$ and define the fields
\[
A(x)=\sqrt{2\log2}U_n(x) \quad\mbox{and} \quad B(x)=\sqrt
{2}Y_{n\log2}(x)+\sqrt{2C\log2}Z.
\]
We have $\mathbb{E}(A(x)A(y))\leq\mathbb{E}(B(x)B(y))$ for all $x,y$.
We then apply the convexity inequality to the fields $A$ and $B$ with
the convex function $G(x)=n^{q({1}/{2})}x^q$ for $q<1$, to get
\begin{eqnarray*}
&& \mathbb{E}\bigl(e^{q\sqrt{2C\log2}Z-qC\log2}\bigr)\mathbb{E}\bigl(S_{n}^q
\bigr)
\\
&&\qquad \leq\mathbb{E} \biggl(n^{q({1}/{2})} \biggl(\int_0^1
e^{ \sqrt
{2\log2}
U_n(x)-\log2\mathbb{E}(U_n(x)^2)}\,\d x \biggr)^q \biggr).
\end{eqnarray*}
Comparing with the notation of \cite{bknsw12}, we see that the quantity
on the right here is simply $\mathbb{E}((n^{{1}/{2} }Z_{1,n})^q)$,
the $q$th moment of the total mass of the correctly renormalized
critical Mandelbrot cascade measure. As noted in \cite{bknsw12}, the
fact that this is uniformly bounded in $n$ for a fixed $q<1$ follows
from \cite{ma11,we11}. Thus, $\mathbb{E}(S_{n}^q)$ is also uniformly
bounded in $n$ for $q<1$. So, recalling that $s_n=\sqrt{\epsilon\log
n}$ and $k_n=\frac{1}{\log2}\log n$, we conclude that for any
$\epsilon
\in(0,1)$, there are constants $C(\frac{\epsilon}{2})$ and
$C(\epsilon
)$ so that if we take $n$ large enough, then $1-\mathbb{P}(\mathcal
{B}_{n})\leq C(\frac{\epsilon}{2})(2Ck_n e^{2\sqrt{2}s_n})^{1-(
{\epsilon}/{2})} n^{(1-({\epsilon}/{2}))(\gamma-({1}/{2}))}\le
C(\epsilon) n^{(1-\epsilon)(\gamma-({1}/{2}))}$. Thus, by (\ref
{eqPmax}) all we are left with is to estimate the Laplace transform of $S_n$.

We make use of the following formula, valid for all nonnegative random
variables $X$:
\[
1-\mathbb{E}\bigl(\exp(-\alpha X)\bigr)=\int_0^\infty
\alpha e^{-\alpha
t}\mathbb{P}(X\geq t)\,\dd t.
\]
In this formula, we set $\alpha=2Ce^{2\sqrt{2}s_n}k_n n^{(\gamma
-({1}/{2}))}$ and $X=S_n$. Recalling from the argument above that $\mathbb
{E}(S_{n}^q)$ is uniformly bounded in $n$ for $q<1$, by Chebyshev's
inequality we see that for any $q<1$
\[
\mathbb{P}(S_{n}\geq t)\leq C_qt^{-q}.
\]
Making the change of variable $\tau=\alpha t$, we get
\[
1-\mathbb{E}\bigl(e^{-\alpha S_{n}}\bigr)\leq C_q
\alpha^{q}\int_0^\infty
e^{-\tau
} \tau^{-q}\,\dd\tau.
\]
Recalling again that $s_n=\sqrt{\epsilon\log n}$ and $k_n=\frac
{1}{\log
2}\log n$, we see that since the integral converges, we can take $q$ so
close to one that we get
\[
1-\mathbb{E}\bigl(e^{-\alpha S_n}\bigr)\leq C'(\epsilon)
n^{(\gamma-({1}/{2}))(1-\epsilon)},
\]
which completes the proof of Lemma \ref{lemmamodulus-probability}.\quad\qed
\end{longlist}\noqed
\end{pf}

Theorem \ref{teomodulus} now follows quickly. We first prove the
analogous statement for the measure $\nu_{\sqrt{2}}$.

%
\begin{teo}\label{teonu-modulus}
For any interval $I\subset[0,1]$ and $\gamma<\frac{1}{2}$, almost surely
%
%
\begin{equation}
\nu_{\sqrt{2}}(I)\leq C(\omega) \bigl(\log\bigl(1+|I|^{-1}\bigr)
\bigr)^{-\gamma},
\end{equation}
where $C(\omega)$ is an almost surely finite random constant.
\end{teo}

\begin{pf}
It is enough to restrict to dyadic subintervals. Pick $\gamma\in
(0,\frac
{1}{2})$. Let $l$ be an integer so that $l(\gamma-\frac{1}{2})<-2$. We
then have by Lemma \ref{lemmamodulus-probability} that
\[
\sum_{k=1}^\infty\mathbb{P}\Bigl(\max
_{\sigma\in\Sigma_{k^l}^{(e/o)}} \nu_{\sqrt{2}}(I_\sigma)\geq
k^{-l\gamma}\Bigr)\leq C\sum_{k=1}^\infty
k^{l(({\gamma-({1}/{2})})/{2})}<\infty.
\]
By Borel--Cantelli,
\[
\max_{\sigma\in\Sigma_{k^l}^{(e/o)}}\nu_{\sqrt{2}}(I_\sigma)\leq C(
\omega) k^{-l\gamma}
\]
for a random (almost surely finite) constant $C(\omega)$. Combining the
estimates for even and odd intervals, we get
\[
\max_{\sigma\in\Sigma_{k^l}} \nu_{\sqrt{2}}(I_\sigma)\leq
C'(\omega) k^{-l\gamma}.
\]

We note that $\max_{\sigma\in\Sigma_n}\nu_{\sqrt{2}}(I_\sigma)$ is
decreasing in $n$ so for $k^l\leq n\leq(k+1)^l$ we have
\[
\max_{\sigma\in\Sigma_n}\nu_{\sqrt{2}}(I_\sigma)\leq\max
_{\sigma\in
\Sigma_{k^l}} \nu_{\sqrt{2}}(I_\sigma)\leq
C'(\omega) k^{-l\gamma
}\leq C'(
\omega)2^{l\gamma}n^{-\gamma},
\]
which is the desired result.
\end{pf}

\begin{pf*}{Proof of Theorem \ref{teomodulus}}
From the definition of $\nu_{\sqrt{2}}$, we note that for any interval
$I \subset[0,1]$
\[
\mu_{\sqrt{2}}(I)\leq e^{\sqrt{2}\max_{x\in[0,1]}X_0(x)-1}\nu_{\sqrt{2}}(I),
\]
where $(X_0(x))_{x\in[0,1]}$ is a Gaussian process with a continuous
covariance kernel. The quantity $e^{\sqrt{2}\max_{x\in[0,1]}X_0(x)}$
is almost surely finite, so Theorem \ref{teonu-modulus} implies the result.
\end{pf*}

\section{On the \texorpdfstring{$\mu_{\sqrt{2}}$}{mu_{sqrt{2}}}-almost everywhere local behavior of \texorpdfstring{$\mu_{\sqrt{2}}$}{mu_{sqrt{2}}}}

We consider the following question: what can be said of the size of
smallest possible sets of full $\mu_{\sqrt{2}}$-measure? This question
is partially answered by Theorem \ref{teohausdorff-gauge}, which is
proven in this section.

Let $f\dvtx  \mathbb{N}\to\mathbb{R}^+$ be an ultimately nonincreasing
function tending to
$0$ at infinity and consider the sets
\[
E_n^f = \bigl\{ x\dvtx  \mu_{\sqrt{2}}
\bigl(I_n(x)\bigr) \leq f(n) \bigr\}.
\]
We will determine a class of functions $f$ for which we have
\[
\sum_n \mu_{\sqrt{2}} \bigl(E_n^f
\bigr) < \infty\qquad\mbox{almost surely.}
\]
For a nontrivial result, it is already enough to consider the
expectation of the series above. We fix a sequence $(\eta_n)_{n\ge1}$
taking values in $(0,1)$ and write
\begin{eqnarray*}
\mu_{\sqrt{2}} \bigl(E_n^f \bigr) &=& \int
_0^1 \mathbf{1}_{\{\mu
_{\sqrt
{2}}(I_n(x)) \leq f(n)\}}
\mu_{\sqrt{2}}(\d x) = \sum_{\sigma\in
\Sigma
_n}
\mu_{\sqrt{2}}(I_\sigma) \mathbf{1}_{\{\mu_{\sqrt
{2}}(I_\sigma) \leq
f(n)\}}
\\
& \le&\sum_{\sigma\in\Sigma_n} \mu_{\sqrt{2}}(I_\sigma)
\biggl(\frac
{f(n)}{\mu_{\sqrt{2}}(I_\sigma)} \biggr)^{\eta_n} = \sum
_{\sigma
\in
\Sigma_n} \mu_{\sqrt{2}}(I_\sigma)^{1-\eta_n}
f(n)^{\eta_n}.
\end{eqnarray*}
Let $\epsilon_n=-\frac{\log(f(n))}{n}$ take the form
$ \gamma\sqrt{\frac{\log(n)}{n}+\alpha\frac{\log
\log
(n)}{n}}$ for $n\ge3$, where $\alpha>0$ and $\gamma>0$ are to be
prescribed. Assume $\eta_n=\lambda\epsilon_n$.

Denoting by $W_n$ the $n$th level lognormal factor
\begin{eqnarray*}
W_n &\stackrel{d} {=}&\exp\bigl( \sqrt{2} X_n -
\mathbb{E}X_n^2 \bigr) \stackrel{d} {=}\exp\biggl(
\sigma_n N - \frac{\sigma_n^2}{2} \biggr),
\\
N &\sim& N(0,1), \qquad\sigma_n^2 = 2n \log2,
\end{eqnarray*}
we have, for each $\sigma\in\Sigma_n$, $\mu_{\sqrt{2}}(I_\sigma
)\stackrel{d}{=}
2^{-n}W_n Y_n$ where $Y_n$ is a copy of $Y$ independent of $W_n$.
Moreover, by Theorem \ref{teotail} we have $\mathbb{E}(Y^{1-\eta
})=O(\eta
^{-1})$ as $\eta\to0^+$. These remarks yield
\begin{eqnarray*}
\mathbb{E}\mu_{\sqrt{2}} \bigl(E_n^f \bigr)
&\le& 2^n 2^{-n(1-\eta
_n)} \mathbb{E} \bigl(W_n^{1-\eta_n}
\bigr) \mathbb{E}\bigl(Y^{1-\eta_n}\bigr)e^{-n\epsilon_n\eta_n}
\\
&\le& C
e^{n (\log
(2) \eta_n^2-\epsilon_n\eta_n) -\log(\eta_n)}.
\end{eqnarray*}
A computation yields for $n\ge3$
\[
n \bigl(\log(2) \eta_n^2-\epsilon_n
\eta_n\bigr) -\log(\eta_n)= \bigl(c+\tfrac{1}{2}
\bigr)\log(n) + \bigl(c \alpha-\tfrac{1}{2} \bigr)\log\log(n) + O(1),
\]
where $c=\log(2) \lambda^2\gamma^2-\lambda\gamma^2$. With the
order of
magnitude chosen for $\epsilon_n$, taking $c=-\frac{3}{2}$ is optimal
in view of making $\sum_{n\ge1} \mathbb{E}\mu_{\sqrt{2}}
(E_n^f )$
convergent. This condition requires the equation $\log(2) \lambda
^2\gamma^2-\lambda\gamma^2+\frac{3}{2}=0$ to have solutions in
$\lambda
$. This imposes $\gamma\ge\sqrt{6\log(2)}$, hence we choose $\gamma
=\sqrt{6\log(2)}$ to minimize $\epsilon_n$. It then turns out that if
$-\frac{3}{2}\alpha-\frac{1}{2}<-1$, that is, $\alpha>\frac{1}{3}$,
then $\sum_{n\ge1} \mathbb{E}\mu_{\sqrt{2}} (E_n^f
)<\infty$.

Theorem \ref{teohausdorff-gauge} follows from the preceding estimates
by an application of the Borel--Cantelli lemma to the measure $\mu
_{\sqrt{2}}$. As an application of Theorem \ref{teohausdorff-gauge} we
present the following simple corollary.

%
\begin{cor}\label{cor0-dimension}
Almost surely, there exists a set of Hausdorff dimension $0$ that has
full $\mu_{\sqrt{2}}$-measure.
\end{cor}

\begin{pf}
Let
\[
E = \bigl\{ x\dvtx  \mu_{\sqrt{2}}\bigl(I_n(x)\bigr) \geq f(n)
\mbox{ for all but finitely many }n \bigr\},
\]
where $f = f_\alpha$ for some $\alpha> \frac{1}{3}$. Since, by Theorem
\ref{teohausdorff-gauge}, $E$ almost surely has full $\mu_{\sqrt
{2}}$-measure, we only need to show that a.s. it has Hausdorff
dimension $0$.

Let $\{ I_\sigma\}_{\sigma\in\Sigma_n^f} $ be the collection of
dyadic subintervals of $[0,1]$ such that $|\sigma| \geq n$ and $\mu
_{\sqrt{2}}(I_\sigma) \geq f(|\sigma|)$. Clearly, for any $n$, $\{
I_\sigma\}_{\sigma\in\Sigma_n^f}$ is a cover of $E$. But for any $s
> 0$ and sufficiently large $n \in\mathbb{N}$ we have $2^{-({s}/{2})|\sigma
|} \leq\mu_{\sqrt{2}}(I_\sigma)$ for all $\sigma\in\Sigma_n^f$, so
\begin{eqnarray*}
\sum_{\sigma\in\Sigma_n^f} |I_\sigma|^s &=&
\sum_{k \geq n} \sum_{\sigma\in\Sigma_n^f, |\sigma| = k}
|I_\sigma|^s = \sum_{k \geq n} \sum
_{\sigma\in\Sigma_n^f, |\sigma| = k} |I_\sigma|^{s/2}
\bigl(2^{-|\sigma|}\bigr)^{s/2}
\\
& \leq& \sum_{k \geq n} 2^{-k({s}/{2})} \sum
_{\sigma\in\Sigma_n^f,
|\sigma|=k} \mu_{\sqrt{2}}(I_\sigma)
\\
& \leq& \mu_{\sqrt{2}}\bigl([0,1]\bigr) \sum_{k \geq n}
2^{-k({s/2})}. 
\end{eqnarray*}
The last expression tends to $0$ as $n \to\infty$. It follows that for
any $s > 0$ the set $E$ has zero Hausdorff $s$-measure, which implies
the claim.
\end{pf}

\section{Higher dimensions}\label{sehigherdimensions}

In this section, we discuss results corresponding to Theorems \ref
{teotail} and \ref{teomodulus} in a higher-dimensional setting, that
is, for multiplicative chaos measures on $\mathbb{R}^d$ for $d \geq
2$, using
similar methods as in the $d=1$ case. We will focus on the $d=2$ case.
We begin by describing the relevant objects and stating the results,
and we will then sketch the minor differences in the proofs. Finally,
we will make a remark on the higher-dimensional cases $d \geq3$.

Formally, a two-dimensional exactly scale invariant lognormal
multiplicative chaos measure may be constructed by exponentiating a
centered Gaussian field $(X(x))_{x \in\mathbb{R}^2}$ with the
covariance $\mathbb{E}
X(x) X(y) = \log^+ \frac{r}{|x-y|}$, with $r>0$. To make a rigorous
construction (see \cite{bajinrhovar12} Section~A.1), one introduces a
Gaussian process $(X_t(x))_{x \in\mathbb{R}^2,t\geq0}$ 
with covariance:
%
%
\begin{eqnarray}
\label{eq2dcov}
&& \mathbb{E}\bigl(X_t(x)X_s(y)\bigr)
\nonumber\\[25pt]\\[-38pt]
&&\qquad =
\cases{ 0, &\quad$|x-y|>r$,
\vspace*{5pt}\cr
\displaystyle\log\frac{r}{|x-y|},
&\quad$re^{-t\wedge s}<|x-y|\leq r$,
\vspace*{5pt}\cr
t\wedge s+2 \biggl(1-\sqrt{
\frac{|x-y|}{r}e^{t\wedge s}} \biggr), &\quad$|x-y|\leq r
e^{-t\wedge s}$.}\nonumber
\end{eqnarray}
It follows from \cite{drsv12-1,drsv12-2} (see Remark 3 in \cite
{drsv12-2} in particular) that a nontrivial critical measure $\mu$
exists (the critical point being $\beta_c=2$) and it can be written as
%
%
\begin{equation}
\label{eq2d-exactly-scale-invariant} \mu(\dd x)=\lim_{t\to\infty}\sqrt
{t}e^{2 X_t(x)-2(t+2)}\,\dd x,
\end{equation}
where the limit is taken weakly in probability. The measure can also be
constructed through the derivative martingale measure. This measure is
exactly scale invariant, that is, for any $\lambda<1$
\[
\bigl(\mu(\lambda A)\bigr)_{A\in\mathcal{B}(B_{{r}/{2}})}\stackrel
{d} {=}\lambda^2
e^{2X_\lambda-2\mathbb{E}(X_\lambda^2)}\bigl(\mu(A)\bigr)_{A\in\mathcal
{B}(B_{{r}/{2}})},
\]
where $B_{{r}/{2}}$ is any disk of radius $\frac{r}{2}$,
$\mathcal
{B}(B_{{r}/{2}})$ denotes its Borel subsets and $X_\lambda$ is a
centered Gaussian with variance $\log\frac{1}{\lambda}$ and as in the
one-dimensional case, it is independent of $(\mu(A))_{A\in\mathcal
{B}(B_{{r}/{2}})}$. The parameter $r$ plays the role of a scale
parameter. We fix $r=1$ from now on.

Our proof of Theorem \ref{teotail} is robust in the sense that in
addition to exact scale invariance, very little extra information on
the exponentiated field $(X_t(x))$ is used. Indeed, we will prove the
following theorem.

%
\begin{teo}\label{teotail2d}
Let $Q=[0,a]^2$ with $a\leq1$ 
and write $Q_1=[0,\frac{a}{2}]^2$. Then
\[
\lim_{\lambda\to\infty}\lambda\mathbb{P}\bigl(\mu(Q)>\lambda\bigr)=c
\]
for
\[
c=\frac{2}{\log2}\mathbb{E} \biggl(\mu(Q_1)\log\frac{\mu(Q)}{\mu
(Q_1)}
\biggr) < \infty.
\]
\end{teo}

\begin{rem}\label{remsquares} Using different values of $a$ and $r$,
we obtain upper and lower bounds of similar form for disks (or any
other compact set containing an open set) instead of squares. Also,
this result can be used to obtain similar bounds for measures other
than the exactly scale invariant one (e.g., by controlling the
Radon--Nikodym derivative).
\end{rem}

For our proof of the modulus of continuity, we needed a further
decorrelation property of the family of fields $(X_t(x))$ and the
$\star
$-scale invariant measure was more convenient than the exactly scale
invariant one. We define a corresponding one in two dimensions:
consider $Y_t(x)=X_t(x)-X_0(x)$. Again from \cite{drsv12-1,drsv12-2},
it follows that
%
%
\begin{equation}
\nu(\dd x)=\lim_{t\to\infty} \sqrt{t}e^{2Y_t(x)-2t}\,\dd x\label{normlimit}
\end{equation}
exists when the limit is taken weakly in probability, that the limit is
nontrivial and that it has the $\star$-scaling property. The $\star
$-scaling property is a consequence of the fact that for $0 < t < t'$,
the field $Y$ may be decomposed as $Y_{t'}(x) = Y_t(x) + Y_{t,t'}(x)$,
where $Y_{t,t'}$ is a scaled copy of $Y_{t'-t}$ that is sampled
independently of $Y_t$. Especially, $Y_{t,t'}(x)$ is also independent
of $Y_{t,t'}(y)$ for $|x-y| \geq e^{-t}$. This decomposition property
was crucial and also sufficient for the proof of Theorem \ref
{teomodulus}, so without further comment have the following theorem.

%
\begin{teo}\label{teomodulus2d}
Let $Q=[0,a]^2$ with $a\leq1$ 
and $\gamma<\frac{1}{2}$. Then
\[
\nu(Q)\leq C(\omega) \bigl(\log\bigl(1+|Q|^{-1}\bigr)
\bigr)^{-\gamma}
\]
for some random, almost surely finite, constant $C(\omega)$.
\end{teo}

Again, the result readily extends to other sets besides squares, and
also to other measures such as $\mu$.

We now sketch how the proof of Theorem \ref{teotail} should be adapted
in order to prove Theorem \ref{teotail2d}.

First, a fundamental part of our proof of Theorem \ref{teotail} was
that we were able to write
\[
Y = \mu\bigl([0,1]\bigr) = W_0 Y_0 + W_1
Y_1,
\]
where for $i=1,2$, $Y_i \stackrel{d}{=}Y$ and
$W_i \stackrel{d}{=}\frac{1}{4} e^{\sqrt{2 \log2} N }$, where $N$
is normal, and
$W_i$ is independent of $Y_i$. This decomposition followed from the
explicit white noise representation of the field $X_t(x)$
(see the \hyperref[app]{Appendix}) which is lacking in dimension two. In the \hyperref[app]{Appendix},
we prove the following replacement.
%

\begin{lemma}\label{lemma2ddecomposition} Let $Y=\mu(Q)$ and
$Q=[0,a]^2=\bigcup_{i=1}^4Q_i$ where $Q_i$ are squares of side $a/2$. By
possibly enlarging the probability space where the process $(X_t(x))_{x
\in\mathbb{R}^2,t\geq0}$ is defined, we
may write
\begin{eqnarray*}
Y&=& \sum_{i=1}^4 \mu(Q_i)=
\sum_{i=1}^4W_i
Y_i,
\end{eqnarray*}
where for each $i$, $Y_i \stackrel{d}{=}Y$, $W_i \stackrel{d}{=}\frac
{1}{16}e^{2\sqrt
{\log2}N}$ with $N$ a standard
normal variable, and $Y_i$ is independent of $W_i$.
\end{lemma}

With this input, adapting Lemma \ref{lemmamixedmoments} to the
higher-dimensional context turns out to be the only significant task.

\begin{pf*}{Proof of Theorem \ref{teotail2d}}
Using Lemma \ref{lemma2ddecomposition}, we may define the Peyri\`ere
measure $\mathbb{Q}$ on $\Omega\times\{1,2,3,4\}$ by setting
\[
\mathbb{E}_\mathbb{Q}f(\omega,j)=\mathbb{E}\sum
_{j=1}^4 W_j(\omega)f(\omega,j),
\]
and then we may define the random variables $\widetilde{Y}(\omega,j) =
Y_j(\omega)$, $\widetilde{W}(\omega,j) = W_j(\omega)$ and
$\widetilde
{B}(\omega,j) = \sum_{i \neq j} W_i(\omega)Y_i(\omega)$. From this
point on the proof of Theorem \ref{teotail} may be followed with only
cosmetic modifications. Lemma \ref{lemmaqb-variables} holds true, the
measure $\nu$ may be defined exactly as in (\ref{eqnu}) and one
obtains the Poisson equation (\ref{eqpoisson}). To apply Proposition
\ref{proppoisson-equation}, we only need to check there is an analogue
of Lemma \ref{lemmamixedmoments} in the two-dimensional setup. Note
that even though Lemma \ref{lemmamixedmoments} holds for all $h \in
(0,1)$, for the tail result to hold it is sufficient to have the result
for $h \in(0,\frac{1}{2}+\varepsilon)$ for some $\varepsilon> 0$.
This is proven next as Lemma \ref{lemmamixedmoments-2d}.
\end{pf*}

%
\begin{lemma}\label{lemmamixedmoments-2d}
For any $h \in(0,\frac{1}{2}+\frac{1}{2\sqrt{2}})$,
\[
\mathbb{E}\bigl(\mu(Q_1)^h \mu(Q \setminus
Q_1)^h\bigr) < \infty.
\]
\end{lemma}

\begin{pf}
The idea of the proof of Lemma \ref{lemmamixedmoments} may be applied,
but some differences arise from the fact that the boundary points
common to both $Q_1$ and $Q \setminus Q_1$ are two line segments rather
than just one point. We start by noting that Lemma~\ref{lemmamixedmoments-disjoint} has an analogue in this setting, with
exactly the same proof: for two Borel sets $A, B \subset\mathbb{R}^2$
separated by a positive distance, we have $\mathbb{E}(\mu(A)^h \mu
(B)^h) <
\infty$ for any $h \in(0,1)$.

By subadditivity, we may estimate
\begin{eqnarray*}
&& \mathbb{E}\bigl(\mu(Q_1)^h \mu(Q \setminus
Q_1)^h\bigr)
\\
&&\qquad \leq\mathbb{E}\bigl(\mu(Q_1)^h
\mu(Q_2)^h\bigr) + \mathbb{E}\bigl(\mu(Q_1)^h
\mu(Q_3)^h\bigr) + \mathbb{E}\bigl(\mu(Q_1)^h
\mu(Q_4)^h\bigr).
\end{eqnarray*}
Suppose that $Q_2$ and $Q_3$ are the squares that share a boundary
segment with $Q_1$. Then the first two terms on the right are equal and
we need to estimate two different types of terms.

Let us first consider $Q_1=[0,\frac{a}{2}]^2 =: P_1$ and $Q_4 = [\frac
{a}{2},a]^2 =: R_1$. We then decompose
\begin{eqnarray*}
P_1 \times R_1 &=& \biggl( \biggl[\frac{a}{4},
\frac{a}{2} \biggr]^2 \times\biggl[\frac{a}{2},
\frac{3a}{4} \biggr]^2 \biggr) \cup A_1
\\
&=:& (P_2 \times R_2) \cup A_1,
\end{eqnarray*}
where $A_1 = (P_1 \times R_1) \setminus(P_2 \times R_2)$. We note that
$P_2\times R_2$ is simply a scaled and translated version of $P_1\times
R_1$, so we can repeat this procedure. We obtain
%
%
\begin{equation}
\label{eqsquare-decomposition} P_1 \times R_1 = \biggl\{\biggl(
\frac{a}{2},\frac{a}{2}\biggr)\biggr\} \cup\bigcup
_{k=1}^\infty A_k, 
\end{equation}
where $P_{k+1}$ is a square of side length $2^{-k-1}a$ with upper right
corner at $(\frac{a}{2},\frac{a}{2})$ and $R_{k+1}$ is a square of side
length $2^{-k-1}a$ with lower left corner at $(\frac{a}{2},\frac
{a}{2})$. Moreover, the $A_i$ are mutually disjoint and disjoint from
$P_{k+1}\times R_{k+1}$, and $A_k$ is a scaled and translated version
of $A_1$ with the scale factor $2^{-k+1}$. The set $A_1$ is a finite
union of products of two sets with positive distance. Using Lemma \ref
{lemmamixedmoments-disjoint}, we see that $\mathbb{E}((\mu\otimes
\mu
)(A_1)^h)<\infty$, and by exact scaling we have
\[
(\mu\otimes\mu) (A_k)\stackrel{d} {=} 2^{4(-k+1)}e^{4 X_k - 4 \mathbb{E}
X_k^2}(
\mu\otimes\mu) (A_1).
\]
Thus, by subadditivity, the decomposition (\ref
{eqsquare-decomposition}) yields
\[
\mathbb{E}\bigl(\mu(P_1)^h \mu(R_1)^h
\bigr) \leq\mathbb{E}\bigl((\mu\otimes\mu) (A_1)^h\bigr)
\sum_{k=1}^\infty2^{-4(k-1)h}
e^{(8h^2-4h)\mathbb{E}(X_k^2)}.
\]
Since $\mathbb{E}X_k^2 = k \log2$, we see that the series converges
for any
$h\in(0,1)$. We also made use of the fact that almost surely $
(\frac{1}{2},\frac{1}{2} )$ is not an atom.

Consider next the case $P_1:= [0,\frac{a}{2}]^2 = Q_1$ and $R_1:=
[\frac{1}{2},1] \times[0,\frac{1}{2}] = Q_2$. We may then write
\[
P_1 \times R_1 = \bigl(P_2^u
\times R_2^u\bigr) \cup\bigl(P_2^u
\times R_2^d\bigr) \cup\bigl(P_2^d
\times Q_2^u\bigr) \cup\bigl(P_2^d
\cup Q_2^d\bigr) \cup A_1,
\]
where $P_2^u$ is the upper half of $[\frac{a}{4},\frac{a}{2}]\times
[0,\frac{a}{2}]$, $P_2^d$ its lower half and similarly for $R$. The set
$A_1$ is what remains, and again it is a finite union of products of
two sets whose distance is positive. The\vspace*{1pt} terms corresponding to
$P_2^u\times R_2^d$ and $P_2^d \times R_2^u$ are of the form we
considered already and the sets $P_2^u \times R_2^u$ and $P_2^d \times
R_2^d$ are scaled and translated copies of $P_1 \times R_1$. We repeat
this decomposition for $P_2^u\times R_2^u$ and $P_2^d\times R_2^d$ and
iterate. At the $k$th iteration, we have $2^k$ sets of the form
$[0,\frac{a}{2}]^2 \times[\frac{a}{2},a]^2$ scaled by $2^{-k}$ and
having pairwise disjoint interiors, and also $2^{k-1}$ copies of $A_1$
with disjoint interiors, scaled by $2^{-k+1}$. Finally, we also have
$2^k$ terms that are scaled and translated copies of $P_1 \times R_1$,
which are then further decomposed in the $k+1$th step. Using exact
scaling, subadditivity and the fact that the $\mu$-mass of the boundary
segments is almost surely zero, we obtain
%
\begin{eqnarray*}
\mathbb{E}\bigl(\mu(P_1)^h\mu(R_1)^h
\bigr) &\leq& C \mathbb{E}\bigl((\mu\otimes\mu) (A_1)^h
\bigr) \sum_{k=1}^\infty2^{k(1-8h+8h^2)}
\\
&&{} + C' \mathbb{E} \bigl(\mu(Q_1)^h
\mu(Q_4)^h \bigr) \sum_{k=1}^\infty
2^{k(1-8h+8h^2)}.
\end{eqnarray*}
The series converge for $\frac{1}{2}-\frac{1}{2\sqrt{2}}<h<\frac
{1}{2}+\frac{1}{2\sqrt{2}}$, completing the proof of the lemma.
\end{pf}

We close this section by commenting on the case
$d \geq3$. 
It is known (\cite{rhovar10}) that exactly scale invariant
multiplicative chaos measures exist in any dimension, but in the known
cases, the associated Gaussian field has long range correlations for
$d\geq3$ (i.e., the covariance does not have compact support) and due
to this the existence of a
nontrivial critical measure is as of yet an open question. This being
said, such correlations played no role in our proof of Theorem \ref
{teotail}. Indeed, if one could establish the limit (\ref{normlimit})
the proof of Theorem \ref{teotail2d} would also extend to the case $d
\geq3$, with only the combinatorics involved in establishing analogues
of Lemma \ref{lemmamixedmoments-2d} getting slightly more cumbersome.


For the modulus of continuity, the long range correlations, and more
specifically the lack of decompositions of the approximating fields
with the required decorrelation properties, are more problematic and
our proof does not work as it is. On the other hand, in any dimension
there exists a $\star$-scale invariant critical measure which does not
have long range correlations. Thus, a possible way to proceed is to try
to prove the corresponding tail result for this measure.

\begin{appendix}
\section*{Appendix: Scale invariance properties}\label{app}

In this section, we give the computations leading to the statements
(\ref{eqx-exact-scale-invariance}) and (\ref
{eqexact-scale-invariance}) on the exact scale invariance of the field
$X$ and of the measure $\mu_{\sqrt{2}}$. We also discuss the $\star
$-scaling relation for the measure $\nu_{\sqrt{2}}$ given in (\ref
{eqstar-scaling}), and finally prove Lemma~\ref{lemma2ddecomposition}.

\subsection{Scaling properties for critical one-dimensional measures}

%
\begin{prop}\label{propexact-scale-invariance}
The random measure $\mu_{\sqrt{2}}$ satisfies the exact scale
invariance property (\ref{eqexact-scale-invariance}), that is, for any
interval $I \subset[0,1]$
\[
\mu_{\sqrt{2}} \lfloor I \stackrel{d} {=}|I| e^{\sqrt{2} X(I) -
\mathbb{E}X(I)^2} \mu
_{\sqrt{2}}^I,
\]
where $\mu_{\sqrt{2}} \lfloor I$ denotes the restriction of $\mu
_{\sqrt
{2}}$ onto $I$ and $\mu_{\sqrt{2}}^I$ is a random measure independent
of $X(I)$ with the law given by
\[
\bigl(\mu_{\sqrt{2}}^I(A) \bigr)_{A \in\mathcal{B}(I)} \stackrel{d}
{=} \bigl(\mu_{\sqrt{2}}\bigl(|I|^{-1}A\bigr)
\bigr)_{A \in\mathcal{B}(I)}.
\]
\end{prop}

\begin{rem*}
Writing the scaling relation simultaneously
for a set $\{I_j\}$ of subintervals of $[0,1]$, one has
\[
\bigl( \mu_{\sqrt{2}} \lfloor I_j \bigr)_j
\stackrel{d} {=} \bigl( |I_j| e^{\sqrt
{2} X(I_j) - \mathbb{E}X(I_j)^2}
\mu_{\sqrt{2}}^{I_j} \bigr)_j,
\]
where the $\mu_{\sqrt{2}}^{I_j}$ are random measures such that for
each $j$,
\[
\bigl( \mu_{\sqrt{2}}^{I_j}(J) \bigr)_{J \in\mathcal{B}(I_j)}
\stackrel{d} {=} \bigl( \mu_{\sqrt{2}}\bigl(|I_j|^{-1}
J\bigr) \bigr)_{J \in\mathcal{B}(I_j)} \quad\mbox{and} \quad\mu_{\sqrt{2}}^{I_j}
\perp\bigl\{ X(A) \bigr\}_{A
\subset
\mathcal{C}(I_j)}.
\]
However, we stress that for subintervals of the unit interval, for $j
\neq k$ the measure $\mu_{\sqrt{2}}^{I_j}$ is not independent either of
$\mu_{\sqrt{2}}^{I_k}$ or $X(I_k)$.
\end{rem*}

\begin{pf*}{Proof of Proposition \ref{propexact-scale-invariance}}
We first show that (\ref{eqx-exact-scale-invariance}) holds. Consider,
for notational convenience, the interval $I = [0,y]$ with $0 < y < 1$.
By definition, for $t \geq\log1/y$ we have
\[
\bigl(X_t(x)\bigr)_{x \in I} = \bigl(X(I) +
X_t^I(x)\bigr)_{x \in I}.
\]
Therefore, it suffices to check that
\[
\bigl(X_t^I(x)\bigr)_{x \in I} \stackrel{d} {=}
\bigl(X_{t-\log1/y}(x/y)\bigr)_{x \in I}
\]
and since the processes are Gaussian, it is enough to consider the
covariance structures. Checking that the covariances of the processes
are the same is demonstrated in Figure~\ref{figcc-selfsimilarity}.

%
\begin{figure}

\includegraphics{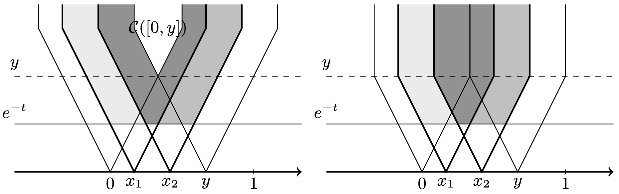}

\caption{\emph{Left.} The sets $\mathcal{C}_t(x_1) \setminus
\mathcal{C}_t(x_2)$
and $\mathcal{C}_t(x_2) \setminus\mathcal{C}_t(x_1)$ are shaded
light gray, while
the intersection $ (\mathcal{C}_t(x_1) \cap\mathcal
{C}_t(x_2) ) \setminus
\mathcal{C}([0,y])$ is dark gray. The law of the Gaussian process
$(X_t^I(x))_{x \in[0,y]}$ is determined by the hyperbolic areas of
these sets for all pairs $(x_1,x_2) \in[0,y]^2$. The set $\mathcal{C}
([0,y])$, contained in every $\mathcal{C}_t(x)$ for $x \in[0,y]$, has been
left white. \emph{Right.} Closing the gap left by the set $\mathcal{C}
([0,y])$ does not affect the hyperbolic areas of any of the shaded
regions. Scaling this picture by $1/y$ also leaves the hyperbolic areas
invariant, giving the distributional equality $(X_t^I(x))_{x \in I}
\stackrel{d}{=}(X_{t-\log1/y}(x/y))_{x \in I}$.}
\label{figcc-selfsimilarity}
\end{figure}

Showing the exact scale invariance of $\mu_{\sqrt{2}}$ is now simple,
as one only needs to note that the measure-defined analogously to the
subcritical measures vanishes: for any intervals $J \subset I \subset
[0,1]$ we have
\begin{eqnarray*}
\mu_{\sqrt{2}}(J) &=& \lim_{t \to\infty} \int
_J \bigl( \sqrt{2}(t+1) - X_t(x) \bigr)
e^{\sqrt{2} X_t(x) - \mathbb{E}X_t(x)^2} \,\d x
\\
&=& \lim_{t \to\infty} \int_J \bigl( \sqrt{2}
\mathbb{E}X(I)^2 - X(I) \bigr) e^{\sqrt{2} X_t(x) - \mathbb{E}X_t(x)^2}
\,\d x
\\
&&{}+ \lim_{t \to\infty} \int_J \bigl(
\sqrt{2}\bigl(t+1-\mathbb{E}X(I)^2\bigr) - X_t^I(x)
\bigr) e^{\sqrt{2} X_t(x) - \mathbb{E}X_t(x)^2} \,\d x
\\
&=& 0
\\
&&{}+ e^{\sqrt{2} X(I) - \mathbb{E}X(I)^2}
\\
&&\quad{}\times  \lim_{t \to\infty} \int
_J \bigl( \sqrt{2}\bigl(t+1-\mathbb{E}X(I)^2
\bigr) - X_t^I(x) \bigr) e^{\sqrt{2}
X_t^I(x) - \mathbb{E}
X_t^I(x)^2} \,\d x
\\
&=:& |I| e^{\sqrt{2} X(I) - \mathbb{E}X(I)^2} \mu^I\bigl(|I|^{-1} J\bigr),
\end{eqnarray*}
where $\mu^I$ a random measure with the law of $\mu$ and independent of
$X(I)$. Note that the measure $\mu^I$ defined here depends on the field
$X$ only through the processes $(X_t^I(x))_{x \in I}$, $t > 0$. This
observation implies the statement on the simultaneous scaling relations
for a set of intervals $\{I_j\}$.
\end{pf*}

We then consider $\star$-scale invariance, as defined in \cite
{ARhVa11}, and the measure $\nu_{\sqrt{2}}$ defined for the proof of
Theorem \ref{teomodulus}. A random measure $\nu$ on $[0,1]$ is called
$\star$-scale invariant on scale $\epsilon\in(0,1]$ if there exist a
process $(\omega_\epsilon(x))_{x \in[0,1]}$ and a random measure
$\nu
^{\epsilon}$ that are independent of each other and satisfy
\[
\bigl( \nu(A) \bigr)_{A \in\mathcal{B}([0,1])} \stackrel{d} {=} \biggl
( \epsilon\int
_A e^{\omega_\epsilon(x)} \,\d\nu^{\epsilon}(x)
\biggr)_{A \in
\mathcal{B}([0,1])}
\]
and
\[
\bigl( \nu^{\epsilon}(A) \bigr)_{A \in\mathcal{B}([0,1])} \stackrel{d}
{=} \bigl( \nu
\bigl(\epsilon^{-1} A\bigr) \bigr)_{A \in\mathcal{B}([0,1])}.
\]

\begin{figure}[t]

\includegraphics{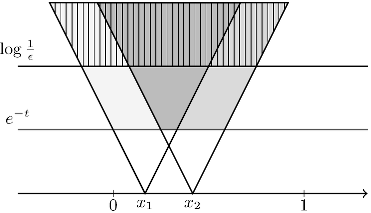}

\caption{The\vspace*{1pt} cones $\mathcal{C}_t(x_1)$ and $\mathcal{C}_t(x_2)$
have been shaded
gray, with the parts in $\mathcal{C}_{\log({1}/{\epsilon})}(x_1)$ and
$\mathcal{C}_{\log({1}/{\epsilon})}(x_2)$ highlighted. By\vspace*{-2pt} scaling
the part
of the picture below the line $\log\frac{1}{\epsilon}$ by $\epsilon
^{-1}$ we get the equality of distributions
$(Y_t(x) - Y_{\log ({1}/{\epsilon})}(x))_{x \in[0,1]} \stackrel{d}{=}
(Y_{t-\log({1}/{\epsilon})}(\epsilon^{-1} x))_{x \in[0,1]}$.
This immediately implies (\protect\ref{eqy-scale-invariance}), since the
process $(Y_{\log ({1}/{\epsilon})}(x))_{x \in[0,1]}$ is independent
of $(Y_t(x) - Y_{\log ({1}/{\epsilon})}(x))_{x \in[0,1]}$.}
\label{figstar-scaling}
\end{figure}

The measure
\[
\nu_{\sqrt{2}}(\d x) = \lim_{t \to\infty} \sqrt{t}
e^{\sqrt{2} Y_t(x)
- \mathbb{E}Y_t(x)^2} \,\d x,
\]
where $Y_t(x) = X_t(x) - X_0(x) = W(\mathcal{C}_t(x) \setminus
\mathcal{C}_0(x))$,
is $\star$-scale invariant on every scale $\epsilon\in(0,1]$ with
\[
\omega_\epsilon(x) = \sqrt{2} Y_{\log({1}/{\epsilon})}(x) + \log
\epsilon.
\]
This can be seen by first deducing the scale invariance property
%
%
\begin{equation}
\label{eqy-scale-invariance} \bigl( Y_t(x) \bigr)_{x \in[0,1]} \stackrel
{d} {=}
\bigl( Y_{\log
({1}/{\epsilon})}(x) + Y_{t - \log({1}/{\epsilon})}'\bigl(
\epsilon^{-1} x\bigr) \bigr)_{x \in[0,1]},
\end{equation}
where $Y'$ is an independent realization of the field $Y$, from
Figure~\ref{figstar-scaling} and then performing a computation
analogous to
the one above for $\mu_{\sqrt{2}}$.

%

\subsection{Joint exact scaling property in two dimensions}
\mbox{}

\begin{pf*}{Proof of Lemma \ref{lemma2ddecomposition}}
For $j=1,\ldots,4$, let $\phi_j\dvtx Q\to Q_j$ be the linear maps that map
the corners of $Q$ to the corners of $Q_j$ by scaling and translating.
We have the following equality in law:
%
%
\begin{equation}
\label{a0} \bigl(X_{t}\bigl(\phi_j(x)\bigr)
\bigr)_{x \in Q,t\geq\log2} \stackrel{d} {=}\bigl( V +X_{t-\log
2}(x)
\bigr)_{x \in Q,t\geq\log2},
\end{equation}
where $V$ is a centered Gaussian variable of variance $\log2$ which is
independent of the process
$(X_{t-\log2}(x))_{x \in Q,t\geq\log2} $. Equation (\ref{a0})
can be readily checked from the form (\ref{eq2dcov}) of the covariance.

We would like to show that by possibly extending our
probability space we can decompose almost surely
%
%
\begin{equation}
\label{a1} X_{t}\bigl(\phi_j(x)\bigr)=
V_j +X^{(j)}_{t-\log2}(x)\qquad\mbox{for }
j=1,2,3,4\mbox{ and } x \in Q,
\end{equation}
where for each $j=1,2,3,4$ the process $(X^{(j)}_{t-\log2}(x))_{x \in
Q,t\geq\log2}$ has the same law as the process
$(X_{t-\log2}(x))_{x \in Q,t \geq\log2}$ and is independent of $V_j$.

As we are interested only in the limit measures, we will need (\ref
{a1}) only for $t$ in some sequence tending to $\infty$. We consider the
countable collection of point evaluations given by
\[
X_{k\log2}(x)\qquad\mbox{where } x\in Q\cap\mathbb{Q}^2,
 k=1,2,3,\ldots
\]
and denote their closed linear span by
\[
\mathcal{H}:=\cspan\bigl(X_{k\log2}(x)\dvtx  x\in Q\cap\mathbb
{Q}^2,  k=1,2,\ldots\bigr).
\]
Thus $\mathcal{H}\subset L^2(\Omega, \mathbb{P})$ is a separable (centered)
Gaussian Hilbert space. By enlarging our probability space, if needed,
we may assume that $(\Omega,\mathcal{F},\mathbb{P})$ supports a centered
Gaussian variable $V$ of variance $\log2$ that is independent of all
elements in $\mathcal{H}$. Set
\[
\mathcal{H}':=\mathcal{H}\oplus
\Span(V).
\]
Consider the closed subspace
\[
\mathcal{G}:=\cspan\bigl(V+X_{k\log2}(x)\dvtx  x\in Q\cap\mathbb
{Q}^2,  k=1,2,\ldots\bigr) \subset\mathcal{H}'.
\]
The dimension of the orthogonal complement of $
\mathcal{G}$ in $\mathcal{H}'$ is either $1$ or zero
since by definition $\cspan(\mathcal{G}\cup\{V\})=\mathcal{H}'$.
Suppose first it is 1 as the latter case is even easier to deal with.
Thus, we may write
\[
\mathcal{H}':=\mathcal{G}\oplus\Span(N),
\]
where $N$ is a centered Gaussian vector of variance $\log2$
independent of all elements in $\mathcal{G}$.

By (\ref{a0}), we have for each $j\in\{ 1,2,3,4\}$ the equality of
joint distributions
%
%
\begin{equation}
\label{eqa2} \bigl(X_{k\log2}\bigl(\phi_j(x)\bigr)
\bigr)_{k\geq1, x\in Q\cap\mathbb{Q}^2} \stackrel{d} {=} \bigl
(V+X_{(k-1)\log2}(x)
\bigr)_{k\geq1, x\in Q\cap\mathbb{Q}^2}.
\end{equation}
This allows us to define linear (not necessarily surjective) isometries
\[
\Psi_j\dvtx \mathcal{H}'=\mathcal{G}\oplus\Span(N)\to
\mathcal{H}'=\mathcal{H}\oplus\Span(V)
\]
as follows. First, set, for $k\geq1$ and $ x\in Q\cap\mathbb{Q}^2$
%
%
\begin{equation}
\label{a11} \Psi_j\bigl(V+X_{(k-1)\log2}(x)
\bigr)=X_{k\log2}\bigl(\phi_j(x)\bigr).
\end{equation}
By (\ref{eqa2}) $\Psi_j$ uniquely extends to an isometry $\Psi
_j\dvtx \mathcal{G}\to
\mathcal{H}'$. Then setting
\[
\Psi_j(N)=V
\]
extends $
\Psi_j$ to the whole of $\mathcal{H}'$.
Note that in case the dimension of the orthogonal complement of $
\mathcal{G}$ in $\mathcal{H}'$ is zero we may omit this last step.

Let us denote
\begin{eqnarray*}
V_j &:=& \Psi_j(V),
\\
X^{(j)}_{k\log2}(x) &:=& \Psi_j
\bigl(X_{k\log2}(x)\bigr)\qquad\mbox{for } k\geq0\mbox{ and }
x\in Q\cap
\mathbb{Q}^2.
\end{eqnarray*}
Since $V$ and $X_{k\log2}(x))$ are independent and $\Psi_j$
is an isometry then $V_j$ is independent of
all the variables $X^{(j)}_{k\log2}(x)$. (\ref{a11}) then gives
%
%
\begin{equation}
\label{eqg2} X_{k\log2}\bigl(\phi_j(x)\bigr) =
V_j + X^{(j)}_{(k-1)\log2}(x)
\end{equation}
for all $k\geq1$, $x\in Q\cap\mathbb{Q}^2$ and $j=1,2,3,4$.

Since the covariance (\ref{eq2dcov}) is H\"older continuous in $x,y$
we may assume that a.s.
$x\to X_{k\log2}(x)$ is continuous.
Since $\Psi_j$ is an isometry the decomposition (\ref{eqg2}) extends
from $x\in Q\cap\mathbb{Q}^2$
to all of $Q$, almost surely.

Consider now, for $k\geq1$, the measures
\[
\mu_k(\dd x):=\sqrt{k\log2} e^{2X_{k\log2}(x)-2\mathbb{E}(X_{k\log
2}(x)^2)}\,\dd x
\]
and for
$k\geq0$
the measures
\[
\mu^{(j)}_k(\dd x):=\sqrt{(k+1)\log2} e^{2X^{(j)}_{k\log
2}(x)-2\mathbb{E}
(X^{(j)}_{k\log2}(x)^2)}\,\dd x.
\]
Using the decomposition
(\ref{eqg2}), we get
%
%
\begin{equation}
\label{a111} \mu_k(Q_j)=\tfrac{1}{4}
e^{2V_j-2\log2} \mu^{(j)}_{k-1}(Q)
\end{equation}
%
and defining
\[
W_j=\tfrac{1}{16}e^{2V_j}
\]
we then get
%
%
\begin{equation}
\label{a1111} \mu_k(Q)=
\sum
_{j=1}^4 W_j \mu^{(j)}_{k-1}(Q).
\end{equation}
Since $
\mu_k\to\mu$ in probability as $k\to\infty$, we infer from
(\ref{a111}) that the variables $\mu^{(j)}_{k-1}(Q) $ converge in probability
to some random variables $Y_j$. Since $\mu^{(j)}_{k-1}(Q) $
has the same distribution as $(\frac{k}{k-1})^{{1}/{2}}\mu
_{k-1}(Q) $, we infer $Y_j \stackrel{d}{=}Y=\mu(Q)$. Hence, taking
limit of
(\ref{a1111}) the desired result follows
as $2V_j$ has variance $4\log2$.
\end{pf*}
\end{appendix}

\section*{Acknowledgements}
We wish to thank the referees for carefully reading the manuscript and
for many helpful suggestions which certainly have improved the quality
of the article.


%

\printaddresses


\begin{thebibliography}{46}

\bibitem{adta07}
%
\begin{bbook}[mr]
\bauthor{\bsnm{Adler},~\bfnm{Robert~J.}\binits{R.~J.}} \AND
\bauthor{\bsnm{Taylor},~\bfnm{Jonathan~E.}\binits{J.~E.}}
(\byear{2007}).
\btitle{Random Fields and Geometry}.
\bpublisher{Springer},
\blocation{New York}.
\bid{mr={2319516}}
\end{bbook}
%
\bptok{imsref}%
\endbibitem

\bibitem{aishi11}
%
\begin{barticle}[mr]
\bauthor{\bsnm{Aidekon},~\bfnm{Elie}\binits{E.}} \AND
\bauthor{\bsnm{Shi},~\bfnm{Zhan}\binits{Z.}}
(\byear{2014}).
\btitle{The {S}eneta--{H}eyde scaling for the branching random walk}.
\bjournal{Ann. Probab.}
\bvolume{42}
\bpages{959--993}.
\bid{doi={10.1214/12-AOP809}, issn={0091-1798}, mr={3189063}}
\bptnote{check year}%
\end{barticle}
%
\bptok{imsref}%
\endbibitem

\bibitem{ARhVa11}
%
\begin{barticle}[mr]
\bauthor{\bsnm{Allez},~\bfnm{Romain}\binits{R.}},
\bauthor{\bsnm{Rhodes},~\bfnm{R{\'e}mi}\binits{R.}} \AND
\bauthor{\bsnm{Vargas},~\bfnm{Vincent}\binits{V.}}
(\byear{2013}).
\btitle{Lognormal {$\star$}-scale invariant random measures}.
\bjournal{Probab. Theory Related Fields}
\bvolume{155}
\bpages{751--788}.
\bid{doi={10.1007/s00440-012-0412-9}, issn={0178-8051}, mr={3034792}}
\bptnote{check year}%
\end{barticle}\vadjust{\goodbreak}
%
\bptok{imsref}%
\endbibitem

\bibitem{AJKS}
%
\begin{barticle}[mr]
\bauthor{\bsnm{Astala},~\bfnm{Kari}\binits{K.}},
\bauthor{\bsnm{Jones},~\bfnm{Peter}\binits{P.}},
\bauthor{\bsnm{Kupiainen},~\bfnm{Antti}\binits{A.}} \AND
\bauthor{\bsnm{Saksman},~\bfnm{Eero}\binits{E.}}
(\byear{2011}).
\btitle{Random conformal weldings}.
\bjournal{Acta Math.}
\bvolume{207}
\bpages{203--254}.
\bid{doi={10.1007/s11511-012-0069-3}, issn={0001-5962}, mr={2892610}}
\end{barticle}
%
\bptok{imsref}%
\endbibitem

\bibitem{BaMu}
%
\begin{barticle}[mr]
\bauthor{\bsnm{Bacry},~\bfnm{E.}\binits{E.}} \AND
\bauthor{\bsnm{Muzy},~\bfnm{J.~F.}\binits{J.~F.}}
(\byear{2003}).
\btitle{Log-infinitely divisible multifractal processes}.
\bjournal{Comm. Math. Phys.}
\bvolume{236}
\bpages{449--475}.
\bid{doi={10.1007/s00220-003-0827-3}, issn={0010-3616}, mr={2021198}}
\end{barticle}
%
\bptok{imsref}%
\endbibitem

\bibitem{BFan05}
%
\begin{barticle}[mr]
\bauthor{\bsnm{Barral},~\bfnm{Julien}\binits{J.}} \AND
\bauthor{\bsnm{Fan},~\bfnm{Ai-Hua}\binits{A.-H.}}
(\byear{2005}).
\btitle{Covering numbers of different points in {D}voretzky covering}.
\bjournal{Bull. Sci. Math.}
\bvolume{129}
\bpages{275--317}.
\bid{doi={10.1016/j.bulsci.2004.05.007}, issn={0007-4497}, mr={2134123}}
\end{barticle}
%
\bptok{imsref}%
\endbibitem

\bibitem{baji12}
\begin{barticle}[mr]
\bauthor{\bsnm{Barral},~\bfnm{Julien}\binits{J.}} \AND
\bauthor{\bsnm{Jin},~\bfnm{Xiong}\binits{X.}}
(\byear{2014}).
\btitle{On exact scaling log-infinitely divisible cascades}.
\bjournal{Probab. Theory Related Fields}
\bvolume{160}
\bpages{521--565}.
\bid{doi={10.1007/s00440-013-0534-8}, issn={0178-8051}, mr={3278915}}
\end{barticle}
%
\bptok{imsref}%
\endbibitem

\bibitem{bajinrhovar12}
%
\begin{barticle}[mr]
\bauthor{\bsnm{Barral},~\bfnm{Julien}\binits{J.}},
\bauthor{\bsnm{Jin},~\bfnm{Xiong}\binits{X.}},
\bauthor{\bsnm{Rhodes},~\bfnm{R{\'e}mi}\binits{R.}} \AND
\bauthor{\bsnm{Vargas},~\bfnm{Vincent}\binits{V.}}
(\byear{2013}).
\btitle{Gaussian multiplicative chaos and {KPZ} duality}.
\bjournal{Comm. Math. Phys.}
\bvolume{323}
\bpages{451--485}.
\bid{doi={10.1007/s00220-013-1769-z}, issn={0010-3616}, mr={3096527}}
\end{barticle}
%
\bptok{imsref}%
\endbibitem

\bibitem{bknsw12}
%
\begin{barticle}[mr]
\bauthor{\bsnm{Barral},~\bfnm{Julien}\binits{J.}},
\bauthor{\bsnm{Kupiainen},~\bfnm{Antti}\binits{A.}},
\bauthor{\bsnm{Nikula},~\bfnm{Miika}\binits{M.}},
\bauthor{\bsnm{Saksman},~\bfnm{Eero}\binits{E.}} \AND
\bauthor{\bsnm{Webb},~\bfnm{Christian}\binits{C.}}
(\byear{2014}).
\btitle{Critical {M}andelbrot cascades}.
\bjournal{Comm. Math. Phys.}
\bvolume{325}
\bpages{685--711}.
\bid{doi={10.1007/s00220-013-1829-4}, issn={0010-3616}, mr={3148099}}
\bptnote{check year}%
\end{barticle}
%
\bptok{imsref}%
\endbibitem

\bibitem{BM1}
%
\begin{barticle}[mr]
\bauthor{\bsnm{Barral},~\bfnm{Julien}\binits{J.}} \AND
\bauthor{\bsnm{Mandelbrot},~\bfnm{Beno{\^{\i}}t~B.}\binits{B.~B.}}
(\byear{2002}).
\btitle{Multifractal products of cylindrical pulses}.
\bjournal{Probab. Theory Related Fields}
\bvolume{124}
\bpages{409--430}.
\bid{doi={10.1007/s004400200220}, issn={0178-8051}, mr={1939653}}
\end{barticle}
%
\bptok{imsref}%
\endbibitem

\bibitem{barhova12}
%
\begin{barticle}[mr]
\bauthor{\bsnm{Barral},~\bfnm{Julien}\binits{J.}},
\bauthor{\bsnm{Rhodes},~\bfnm{R{\'e}mi}\binits{R.}} \AND
\bauthor{\bsnm{Vargas},~\bfnm{Vincent}\binits{V.}}
(\byear{2012}).
\btitle{Limiting laws of supercritical branching random walks}.
\bjournal{C. R. Math. Acad. Sci. Paris}
\bvolume{350}
\bpages{535--538}.
\bid{doi={10.1016/j.crma.2012.05.013}, issn={1631-073X}, mr={2929063}}
\end{barticle}
%
\bptok{imsref}%
\endbibitem

\bibitem{BS09}
%
\begin{barticle}[mr]
\bauthor{\bsnm{Benjamini},~\bfnm{Itai}\binits{I.}} \AND
\bauthor{\bsnm{Schramm},~\bfnm{Oded}\binits{O.}}
(\byear{2009}).
\btitle{K{PZ} in one dimensional random geometry of multiplicative cascades}.
\bjournal{Comm. Math. Phys.}
\bvolume{289}
\bpages{653--662}.
\bid{doi={10.1007/s00220-009-0752-1}, issn={0010-3616}, mr={2506765}}
\end{barticle}
%
\bptok{imsref}%
\endbibitem

\bibitem{bu07}
%
\begin{barticle}[mr]
\bauthor{\bsnm{Buraczewski},~\bfnm{Dariusz}\binits{D.}}
(\byear{2007}).
\btitle{On invariant measures of stochastic recursions in a critical case}.
\bjournal{Ann. Appl. Probab.}
\bvolume{17}
\bpages{1245--1272}.
\bid{doi={10.1214/105051607000000140}, issn={1050-5164}, mr={2344306}}
\end{barticle}
%
\bptok{imsref}%
\endbibitem

\bibitem{bu09}
%
\begin{barticle}[mr]
\bauthor{\bsnm{Buraczewski},~\bfnm{Dariusz}\binits{D.}}
(\byear{2009}).
\btitle{On tails of fixed points of the smoothing transform in the
boundary case}.
\bjournal{Stochastic Process. Appl.}
\bvolume{119}
\bpages{3955--3961}.
\bid{doi={10.1016/j.spa.2009.09.005}, issn={0304-4149}, mr={2552312}}
\end{barticle}
%
\bptok{imsref}%
\endbibitem

\bibitem{cld}
%
\begin{barticle}[auto:STB|2014/02/12|14:17:21]
\bauthor{\bsnm{Carpentier},~\bfnm{D.}\binits{D.}} \AND
\bauthor{\bsnm{Le Doussal},~\bfnm{P.}\binits{P.}}
(\byear{2001}).
\btitle{Glass transition of a particle in a random potential, front
selection in nonlinear RG and entropic phenomena in Liouville and
Sinh--Gordon models}.
\bjournal{Phys. Rev. E (3)}
\bvolume{63}
\bpages{026110}.
\end{barticle}
%
\bptok{imsref}%
\endbibitem

\bibitem{Duphouches}
%
\begin{bincollection}[mr]
\bauthor{\bsnm{Duplantier},~\bfnm{B.}\binits{B.}}
(\byear{2010}).
\btitle{A rigorous perspective on {L}iouville quantum gravity and the
{KPZ} relation}.
In \bbooktitle{Exact Methods in Low-dimensional Statistical Physics and
Quantum Computing}
(\beditor{\bfnm{J.}\binits{J.}~\bsnm{Jacobsen}},
\beditor{\bfnm{S.}\binits{S.}~\bsnm{Ouvry}},
\beditor{\bfnm{V.}\binits{V.}~\bsnm{Pasquier}},
\beditor{\bfnm{D.}\binits{D.}~\bsnm{Serban}}
\AND
\beditor{\bfnm{L.~F.}\binits{L.~F.}~\bsnm{Cugliandolo}}, eds.).
\bseries{Lecture Notes of the Les Houches Summer School}
\bvolume{89}
\bpages{529--561}.
\bpublisher{Oxford Univ. Press},
\blocation{Oxford}.
\bid{mr={2668656}}
\end{bincollection}
%
\bptok{imsref}%
\endbibitem

\bibitem{drsv12-1}
%
\begin{barticle}[auto:STB|2014/02/12|14:17:21]
\bauthor{\bsnm{Duplantier},~\bfnm{B.}\binits{B.}},
\bauthor{\bsnm{Rhodes},~\bfnm{R.}\binits{R.}},
\bauthor{\bsnm{Sheffield},~\bfnm{S.}\binits{S.}} \AND
\bauthor{\bsnm{Vargas},~\bfnm{V.}\binits{V.}}
(\byear{2014}).
\btitle{Critical Gaussian multiplicative chaos: Convergence of the derivative martingale}.
\bjournal{Ann. Probab.}
\bvolume{42}
\bpages{1769--1808}.
\end{barticle}
%
\bptok{imsref}%
\endbibitem

\bibitem{drsv12-2}
\begin{barticle}[mr]
\bauthor{\bsnm{Duplantier},~\bfnm{Bertrand}\binits{B.}},
\bauthor{\bsnm{Rhodes},~\bfnm{R{\'e}mi}\binits{R.}},
\bauthor{\bsnm{Sheffield},~\bfnm{Scott}\binits{S.}} \AND
\bauthor{\bsnm{Vargas},~\bfnm{Vincent}\binits{V.}}
(\byear{2014}).
\btitle{Renormalization of critical {G}aussian multiplicative chaos and {KPZ} relation}.
\bjournal{Comm. Math. Phys.}
\bvolume{330}
\bpages{283--330}.
\bid{doi={10.1007/s00220-014-2000-6}, issn={0010-3616}, mr={3215583}}
\end{barticle}
%
\bptok{imsref}%
\endbibitem

\bibitem{PRL}
%
\begin{barticle}[mr]
\bauthor{\bsnm{Duplantier},~\bfnm{Bertrand}\binits{B.}} \AND
\bauthor{\bsnm{Sheffield},~\bfnm{Scott}\binits{S.}}
(\byear{2009}).
\btitle{Duality and the {K}nizhnik--{P}olyakov--{Z}amolodchikov
relation in {L}iouville quantum gravity}.
\bjournal{Phys. Rev. Lett.}
\bvolume{102}
\bpages{150603, 4}.
\bid{doi={10.1103/PhysRevLett.102.150603}, issn={0031-9007}, mr={2501276}}
\end{barticle}
%
\bptok{imsref}%
\endbibitem

\bibitem{DS08}
%
\begin{barticle}[mr]
\bauthor{\bsnm{Duplantier},~\bfnm{Bertrand}\binits{B.}} \AND
\bauthor{\bsnm{Sheffield},~\bfnm{Scott}\binits{S.}}
(\byear{2011}).
\btitle{Liouville quantum gravity and {KPZ}}.
\bjournal{Invent. Math.}
\bvolume{185}
\bpages{333--393}.
\bid{doi={10.1007/s00222-010-0308-1}, issn={0020-9910}, mr={2819163}}
\end{barticle}
%
\bptok{imsref}%
\endbibitem

\bibitem{duli83}
%
\begin{barticle}[mr]
\bauthor{\bsnm{Durrett},~\bfnm{Richard}\binits{R.}} \AND
\bauthor{\bsnm{Liggett},~\bfnm{Thomas~M.}\binits{T.~M.}}
(\byear{1983}).
\btitle{Fixed points of the smoothing transformation}.
\bjournal{Z.~Wahrsch. Verw. Gebiete}
\bvolume{64}
\bpages{275--301}.
\bid{doi={10.1007/BF00532962}, issn={0044-3719}, mr={0716487}}
\end{barticle}
%
\bptok{imsref}%
\endbibitem

\bibitem{Fan04}
%
\begin{barticle}[mr]
\bauthor{\bsnm{Fan},~\bfnm{Aihua}\binits{A.}}
(\byear{2004}).
\btitle{Limsup deviations on trees}.
\bjournal{Anal. Theory Appl.}
\bvolume{20}
\bpages{113--148}.
\bid{doi={10.1007/BF02901437}, issn={1672-4070}, mr={2095456}}
\end{barticle}\vadjust{\goodbreak}
%
\bptok{imsref}%
\endbibitem

\bibitem{Fan96}
%
\begin{barticle}[mr]
\bauthor{\bsnm{Fan},~\bfnm{Ai~Hua}\binits{A.~H.}}
(\byear{1997}).
\btitle{Sur les chaos de {L}\'evy stables d'indice {$0< \alpha< 1$}}.
\bjournal{Ann. Sci. Math. Qu\'{e}bec}
\bvolume{21}
\bpages{53--66}.
\bid{issn={0707-9109}, mr={1457064}}
\end{barticle}
%
\bptok{imsref}%
\endbibitem

\bibitem{fb}
%
\begin{barticle}[mr]
\bauthor{\bsnm{Fyodorov},~\bfnm{Yan~V.}\binits{Y.~V.}} \AND
\bauthor{\bsnm{Bouchaud},~\bfnm{Jean-Philippe}\binits{J.-P.}}
(\byear{2008}).
\btitle{Freezing and extreme-value statistics in a random energy model
with logarithmically correlated potential}.
\bjournal{J. Phys. A}
\bvolume{41}
\bpages{372001, 12}.
\bid{doi={10.1088/1751-8113/41/37/372001}, issn={1751-8113}, mr={2430565}}
\end{barticle}
%
\bptok{imsref}%
\endbibitem

\bibitem{gu90}
%
\begin{barticle}[mr]
\bauthor{\bsnm{Guivarc'h},~\bfnm{Yves}\binits{Y.}}
(\byear{1990}).
\btitle{Sur une extension de la notion de loi semi-stable}.
\bjournal{Ann. Inst. Henri Poincar\'e Probab. Stat.}
\bvolume{26}
\bpages{261--285}.
\bid{issn={0246-0203}, mr={1063751}}
\end{barticle}
%
\bptok{imsref}%
\endbibitem

\bibitem{K2}
%
\begin{barticle}[mr]
\bauthor{\bsnm{Kahane},~\bfnm{Jean-Pierre}\binits{J.-P.}}
(\byear{1985}).
\btitle{Sur le chaos multiplicatif}.
\bjournal{Ann. Sci. Math. Qu\'{e}bec}
\bvolume{9}
\bpages{105--150}.
\bid{issn={0707-9109}, mr={0829798}}
\end{barticle}
%
\bptok{imsref}%
\endbibitem

\bibitem{K4}
%
\begin{bbook}[mr]
\bauthor{\bsnm{Kahane},~\bfnm{Jean-Pierre}\binits{J.-P.}}
(\byear{1985}).
\btitle{Some Random Series of Functions},
\bedition{2nd} ed.
\bpublisher{Cambridge Univ. Press},
\blocation{Cambridge}.
\bid{mr={0833073}}
\end{bbook}
%
\bptok{imsref}%
\endbibitem

\bibitem{K1}
%
\begin{barticle}[mr]
\bauthor{\bsnm{Kahane},~\bfnm{Jean-Pierre}\binits{J.-P.}}
(\byear{1987}).
\btitle{Positive martingales and random measures}.
\bjournal{Chinese Ann. Math. Ser. B}
\bvolume{8}
\bpages{1--12}.
\bid{issn={0252-9599}, mr={0886744}}
\end{barticle}
%
\bptok{imsref}%
\endbibitem

\bibitem{K3}
%
\begin{barticle}[mr]
\bauthor{\bsnm{Kahane},~\bfnm{Jean-Pierre}\binits{J.-P.}}
(\byear{1987}).
\btitle{Multiplications al\'eatoires et dimensions de {H}ausdorff}.
\bjournal{Ann. Inst. Henri Poincar\'e Probab. Stat.}
\bvolume{23}
\bpages{289--296}.
\bid{issn={0246-0203}, mr={0898497}}
\end{barticle}
%
\bptok{imsref}%
\endbibitem

\bibitem{K5}
%
\begin{barticle}[mr]
\bauthor{\bsnm{Kahane},~\bfnm{Jean-Pierre}\binits{J.-P.}}
(\byear{1990}).
\btitle{Recouvrements al\'eatoires et th\'eorie du potentiel}.
\bjournal{Colloq. Math.}
\bvolume{60/61}
\bpages{387--411}.
\bid{issn={0010-1354}, mr={1096386}}
\end{barticle}
%
\bptok{imsref}%
\endbibitem

\bibitem{kape76}
%
\begin{barticle}[mr]
\bauthor{\bsnm{Kahane},~\bfnm{J.-P.}\binits{J.-P.}} \AND
\bauthor{\bsnm{Peyri{\`e}re},~\bfnm{J.}\binits{J.}}
(\byear{1976}).
\btitle{Sur certaines martingales de {B}enoit {M}andelbrot}.
\bjournal{Adv. Math.}
\bvolume{22}
\bpages{131--145}.
\bid{issn={0001-8708}, mr={0431355}}
\end{barticle}
%
\bptok{imsref}%
\endbibitem

\bibitem{li01}
%
\begin{barticle}[mr]
\bauthor{\bsnm{Liu},~\bfnm{Quansheng}\binits{Q.}}
(\byear{2000}).
\btitle{On generalized multiplicative cascades}.
\bjournal{Stochastic Process. Appl.}
\bvolume{86}
\bpages{263--286}.
\bid{doi={10.1016/S0304-4149(99)00097-6}, issn={0304-4149}, mr={1741808}}
\end{barticle}
%
\bptok{imsref}%
\endbibitem

\bibitem{ma11}
%
\begin{bmisc}[auto:STB|2014/02/12|14:17:21]
\bauthor{\bsnm{Madaule},~\bfnm{T.}\binits{T.}}
(\byear{2011}).
\bhowpublished{Convergence in law for the branching random walk seen
from its tip.
Available at \arxivurl{arXiv:1107.2543}.}
\end{bmisc}
%
\bptok{imsref}%
\endbibitem

\bibitem{marhva13}
%
\begin{bmisc}[auto:STB|2014/02/12|14:17:21]
\bauthor{\bsnm{Madaule},~\bfnm{T.}\binits{T.}},
\bauthor{\bsnm{Rhodes},~\bfnm{R.}\binits{R.}} \AND
\bauthor{\bsnm{Vargas},~\bfnm{V.}\binits{V.}}
(\byear{2013}).
\bhowpublished{Glassy phase and freezing of log-correlated Gaussian potentials.
Available at \arxivurl{arXiv:1310.5574}.}
\end{bmisc}
%
\bptok{imsref}%
\endbibitem

\bibitem{M2}
%
\begin{barticle}[mr]
\bauthor{\bsnm{Mandelbrot},~\bfnm{Benoit}\binits{B.}}
(\byear{1974}).
\btitle{Multiplications al\'eatoires it\'er\'ees et distributions
invariantes par moyenne pond\'er\'ee al\'eatoire}.
\bjournal{C. R. Acad. Sci. Paris S\'er. A}
\bvolume{278}
\bpages{289--292}.
\bid{mr={0431351}}
\end{barticle}
%
\bptok{imsref}%
\endbibitem

\bibitem{M1}
%
\begin{bincollection}[auto:STB|2014/02/12|14:17:21]
\bauthor{\bsnm{Mandelbrot},~\bfnm{B.~B.}\binits{B.~B.}}
(\byear{1972}).
\btitle{Possible refinement of the lognormal hypothesis concerning the
distribution of energy in intermittent turbulence}.
In \bbooktitle{Statistical Models and Turbulence}
(\beditor{\bfnm{M.}\binits{M.}~\bsnm{Rosenblatt}} \and
\beditor{\bfnm{C.~V.}\binits{C.~V.}~\bsnm{Atta}}, eds.).
\bseries{Lectures Notes in Physics}
\bvolume{12}
\bpages{333--351}.
\bpublisher{Springer}, \blocation{New York}.
\end{bincollection}
%
\bptok{imsref}%
\endbibitem

\bibitem{M3}
%
\begin{barticle}[auto:STB|2014/02/12|14:17:21]
\bauthor{\bsnm{Mandelbrot},~\bfnm{B.~B.}\binits{B.~B.}}
(\byear{1974}).
\btitle{Intermittent turbulence in self-similar cascades, divergence of
high moments and dimension of the carrier}.
\bjournal{J. Fluid Mech.}
\bvolume{62}
\bpages{331--358}.
\end{barticle}
%
\bptok{imsref}%
\endbibitem

\bibitem{M89}
%
\begin{bincollection}[mr]
\bauthor{\bsnm{Mandelbrot},~\bfnm{Benoit~B.}\binits{B.~B.}}
(\byear{1989}).
\btitle{Multifractal measures, especially for the geophysicist}.
In \bbooktitle{Fractals in Geophysics}
\bpages{5--42}.
\bpublisher{Birkh\"auser},
\blocation{Basel}.
\bid{mr={1106479}}
\end{bincollection}
%
\bptok{imsref}%
\endbibitem

\bibitem{M97}
%
\begin{bbook}[mr]
\bauthor{\bsnm{Mandelbrot},~\bfnm{Benoit~B.}\binits{B.~B.}}
(\byear{1997}).
\btitle{Fractals and Scaling in Finance: Discontinuity, Concentration, Risk}.
\bpublisher{Springer},
\blocation{New York}.
\bid{doi={10.1007/978-1-4757-2763-0}, mr={1475217}}
\end{bbook}
%
\bptok{imsref}%
\endbibitem

\bibitem{post69}
%
\begin{barticle}[mr]
\bauthor{\bsnm{Port},~\bfnm{Sidney~C.}\binits{S.~C.}} \AND
\bauthor{\bsnm{Stone},~\bfnm{Charles~J.}\binits{C.~J.}}
(\byear{1969}).
\btitle{Potential theory of random walks on {A}belian groups}.
\bjournal{Acta Math.}
\bvolume{122}
\bpages{19--114}.
\bid{issn={0001-5962}, mr={0261706}}
\end{barticle}
%
\bptok{imsref}%
\endbibitem


\bibitem{rhovar10}
%
\begin{barticle}[mr]
\bauthor{\bsnm{Rhodes},~\bfnm{R{\'e}mi}\binits{R.}} \AND
\bauthor{\bsnm{Vargas},~\bfnm{Vincent}\binits{V.}}
(\byear{2010}).
\btitle{Multidimensional multifractal random measures}.
\bjournal{Electron. J. Probab.}
\bvolume{15}
\bpages{241--258}.
\bid{doi={10.1214/EJP.v15-746}, issn={1083-6489}, mr={2609587}}
\end{barticle}
%
\bptok{imsref}%
\endbibitem

\bibitem{rhovar08}
%
\begin{barticle}[mr]
\bauthor{\bsnm{Rhodes},~\bfnm{R{\'e}mi}\binits{R.}} \AND
\bauthor{\bsnm{Vargas},~\bfnm{Vincent}\binits{V.}}
(\byear{2011}).
\btitle{K{PZ} formula for log-infinitely divisible multifractal random
measures}.
\bjournal{ESAIM Probab. Stat.}
\bvolume{15}
\bpages{358--371}.
\bid{doi={10.1051/ps/2010007}, issn={1292-8100}, mr={2870520}}
\end{barticle}
%
\bptok{imsref}%
\endbibitem

\bibitem{RoVa1}
%
\begin{barticle}[mr]
\bauthor{\bsnm{Robert},~\bfnm{Raoul}\binits{R.}} \AND
\bauthor{\bsnm{Vargas},~\bfnm{Vincent}\binits{V.}}
(\byear{2008}).
\btitle{Hydrodynamic turbulence and intermittent random fields}.
\bjournal{Comm. Math. Phys.}
\bvolume{284}
\bpages{649--673}.
\bid{doi={10.1007/s00220-008-0642-y}, issn={0010-3616}, mr={2452591}}
\end{barticle}
%
\bptok{imsref}%
\endbibitem

\bibitem{RoVa2}
%
\begin{barticle}[mr]
\bauthor{\bsnm{Robert},~\bfnm{Raoul}\binits{R.}} \AND
\bauthor{\bsnm{Vargas},~\bfnm{Vincent}\binits{V.}}
(\byear{2010}).
\btitle{Gaussian multiplicative chaos revisited}.
\bjournal{Ann. Probab.}
\bvolume{38}
\bpages{605--631}.
\bid{doi={10.1214/09-AOP490}, issn={0091-1798}, mr={2642887}}
\end{barticle}
%
\bptok{imsref}%
\endbibitem

\bibitem{we11}
%
\begin{barticle}[mr]
\bauthor{\bsnm{Webb},~\bfnm{Christian}\binits{C.}}
(\byear{2011}).
\btitle{Exact asymptotics of the freezing transition of a
logarithmically correlated random energy model}.
\bjournal{J. Stat. Phys.}
\bvolume{145}
\bpages{1595--1619}.
\bid{doi={10.1007/s10955-011-0359-8}, issn={0022-4715}, mr={2863721}}
\end{barticle}
%
\bptok{imsref}%
\endbibitem

\end{thebibliography}
\end{document}